\documentclass[aos,seceqn,nameyear,dvips]{arximspdf}
\usepackage{mathrsfs,graphicx}

\doi{10.1214/10-AOS816}
\volume{38}
\issue{5}
\pubyear{2010}
\firstpage{3028}
\lastpage{3062}

\makeatletter

\newtheorem{theorem}{Theorem}[section]
\newtheorem{pp}[theorem]{Proposition}
\newtheorem{lm}{Lemma}

\newproclaim{remark}{Remark}

\renewcommand{\mathcal}{\mathscr}

\newcommand{\var}{\operatorname{Var}}
\newcommand{\cov}{\operatorname{Cov}}
\newcommand{\tr}{\operatorname{tr}}
\newcommand{\diag}{\operatorname{diag}}
\newcommand{\cid}{\stackrel{d}{\longrightarrow}}
\newcommand{\eid}{\stackrel{d}{=}}
\newcommand{\cip}{\stackrel{p}{\longrightarrow}}
\newcommand{\bdb}{\mathbf{b}}
\newcommand{\bdc}{\mathbf{c}}
\newcommand{\bde}{\mathbf{e}}
\newcommand{\bdg}{\mathbf{g}}
\newcommand{\bdt}{\mathbf{t}}
\newcommand{\bdq}{\mathbf{q}}
\newcommand{\bdx}{\mathbf{x}}
\newcommand{\bdu}{\mathbf{u}}
\newcommand{\bdv}{\mathbf{v}}
\newcommand{\bdw}{\mathbf{w}}
\newcommand{\bdz}{\mathbf{z}}
\newcommand{\CB}{\mathcal{B}}
\newcommand{\CD}{\mathcal{D}}
\newcommand{\CF}{\mathcal{F}}
\newcommand{\CH}{\mathcal{H}}
\newcommand{\CI}{\mathcal{I}}
\newcommand{\CJ}{\mathcal{J}}
\newcommand{\CP}{\mathcal{P}}
\newcommand{\CQ}{\mathcal{Q}}
\newcommand{\CR}{\mathcal{R}}
\newcommand{\CS}{\mathcal{S}}
\newcommand{\CT}{\mathcal{T}}
\newcommand{\CU}{\mathcal{U}}
\newcommand{\CV}{\mathcal{V}}
\newcommand{\CX}{\mathcal{X}}
\newcommand{\CZ}{\mathcal{Z}}
\newcommand{\bbR}{{\mathbb R}}
\newcommand{\E}{{\mathbb E}}
\newcommand{\bbP}{{\mathbb P}}
\newcommand{\BW}{\mathbf{W}}
\newcommand{\BX}{\mathbf{X}}
\newcommand{\BU}{\mathbf{U}}
\newcommand{\bdeta}{{\bolds\eta}}
\newcommand{\bdphi}{{\bolds\phi}}
\newcommand{\bdvartheta}{{\bolds\vartheta}}
\newcommand{\Normal}{\operatorname{Normal}}
\newcommand{\tran}{^{\mathrm{T}}}

\makeatother

\begin{document}
\begin{frontmatter}

\title{Deciding the dimension of effective dimension reduction space
for functional and high-dimensional data}
\runtitle{Deciding the dimension of effective dimension reduction space}

\begin{aug}
\author[A]{\fnms{Yehua} \snm{Li}\thanksref{t1}\ead[label=e1]{yehuali@uga.edu}} and
\author[B]{\fnms{Tailen} \snm{Hsing}\corref{}\thanksref{t2}\ead[label=e2]{thsing@umich.edu}}
\runauthor{Y. Li and T. Hsing}
\affiliation{University of Georgia and University of Michigan}
\address[A]{Department of Statistics\\
University of Georgia\\
Athens, Georgia 30602-7952\\
USA\\
\printead{e1}}
\address[B]{Department of Statistics\\
University of Michigan\\
Ann Arbor, Michigan 48109-1107\\
USA\\
\printead{e2}}
\end{aug}

\thankstext{t1}{Supported by NSF Grant DMS-08-06131.}
\thankstext{t2}{Supported by NSF Grants DMS-08-08993 and DMS-08-06098.}

\received{\smonth{12} \syear{2009}}

%
\begin{abstract}
In this paper, we consider regression models with a Hilbert-space-valued
predictor and a scalar response, where the response depends on
the predictor only through a finite number of projections. The linear
subspace spanned by these projections is called the effective dimension
reduction (EDR) space. To determine the dimensionality of the EDR
space, we focus on the leading\vspace*{2pt} principal component scores
of the
predictor, and propose\vspace*{1pt} two sequential $\chi^2$ testing procedures under
the assumption that the predictor has an elliptically contoured
distribution. We further extend these procedures and introduce a test
that simultaneously takes into account a large number of principal
component scores. The proposed procedures are supported by theory,
validated by simulation studies, and illustrated by a real-data
example. Our methods and theory are applicable to functional data and
high-dimensional multivariate data.
\end{abstract}

%
\begin{keyword}[class=AMS]
\kwd[Primary ]{62J05}
\kwd[; secondary ]{62G20}
\kwd{62M20}.
\end{keyword}
\begin{keyword}
\kwd{Adaptive Neyman test}
\kwd{dimension reduction}
\kwd{elliptically contoured distribution}
\kwd{functional data analysis}
\kwd{principal components}.
\end{keyword}

\end{frontmatter}

\section{Introduction}\label{sec1}
Li (\citeyear{r22}) considered a regression model in which a scalar response
depends on a
multivariate predictor through an unknown number of linear projections, where
the linear space spanned by the directions of the projections was named the
effective dimension reduction (EDR) space of the model. \citet{r22}
introduced a $\chi^{2}$ test
to determine the dimension of the EDR space, and an estimation procedure,
sliced inverse regression (SIR), to estimate the EDR space.
Li's results focused on the case where $p$, the dimension of the
predictor, is
much smaller than~$n$, the sample size. It is not obvious how to extend his
results to high-dimensional multivariate data where $p$ is comparable
to or larger than $n$;
see Remark 5.4 in \citet{r22}.

Regression problems for functional data have drawn a lot of attention recently.
In particular, regression models in which the predictor is functional
while the response is scalar have been extensively investigated;
for linear models, see \citet{r5}, \citet{r28},
\citet{r4} and \citet{r25}; for nonlinear models, see
Hastie and Tibshirani (\citeyear{r18}),
\citet{r6}, \citet{r21} and M\"uller and
Stadtm\"uller (\citeyear{r26}).
Ferr\'e and Yao (\citeyear{r13}, \citeyear{r14}) extended SIR to a
functional-data setting,
and showed
that the EDR space can be consistently estimated under
regularity conditions provided that the true dimension of the space is known;
see also \citet{r14b} and \citet{r14a}.
However, deciding the dimensionality of the EDR space is much more
challenging in that case,
and there has not been a formal procedure to date.

In this paper, we address the problem of deciding the dimensionality
of the EDR space for both functional and high-dimensional
multivariate data. As in \citet{r13}, we adopt the framework
where the predictor takes value in an arbitrary Hilbert space.
To better control the sample information in the (high-dimensional) predictor,
we focus on the sample principal component
scores rather than the raw data. Since the leading principal component
scores optimally explain the
variability in the predictor, it is natural to expect that the leading sample
principal component scores also offer the most relevant information for
the inference problem.
Two statistical tests will be developed for testing whether the
dimension of
the EDR space is larger than a prescribed value;
an estimator of the dimension of the EDR space will then be obtained by
sequentially performing the tests developed.
We will assume that the Hilbert-space-valued predictor has an
elliptically contoured distribution,
a common assumption for inverse regression problems
[cf. \citet{r8}, \citet{r30} and \citet{r13}].
These tests will be
first developed by focusing on a fixed
number of principal component scores; it will be shown that the null
distributions of
the test statistics are asymptotically $\chi^2$.
To address high and infinite-dimensional data, we propose
an ``adaptive Neyman'' test, which combines the information
in a sequence of $\chi^2$ tests corresponding to an increasing number of
principal component scores.


We introduce the background and notation in
Section \ref{sec:fsir}. The main theoretical results and test/estimation
procedures are described in Section
\ref{sec:tests}. Simulation studies are presented in Section
\ref{sec:simulation}, and a real application on near-infrared spectrum
data is presented in Section \ref{sec:data}.
Finally, all of the proofs are collected in \hyperref[app]{Appendix}.

%
\section{Model assumptions and preliminaries}
\label{sec:fsir}

Let $\{X(t),t\in\CI\}$ be a real-valued stochastic
process with an index set $\CI$. Assume that $\bbP(X\in\CH)=1$, where
$\CH$ is some Hilbert space containing functions on $\CI$ and
equipped with inner product
${\langle}\cdot,\cdot{\rangle}$. We do not place any restriction on
$\CI$ and
so $X$ can be
extremely general. For instance, for multivariate data $\CI$ is a
finite set, with
$p$ elements, say, and $\CH$ can then be taken as $\bbR^p$ equipped
with the
usual dot product;
in functional data analysis, $\CH$ is commonly
assumed to be $L^{2}(\CI)$ for some bounded interval $\CI$, with
inner product ${\langle}g,h{\rangle}= \int_\CI g(t)h(t)\,dt$.


Consider the following multiple-index model:
%
%
\begin{equation}\label{eq:fullmodel}
Y=f({\langle}\beta_1, X{\rangle},\ldots, {\langle}\beta_K,
X{\rangle},\varepsilon),
\end{equation}
where $Y$ is scalar, $f$ is an arbitrary function,
$\beta_1,\ldots,\beta_{K}$ are linearly independent elements in $\CH$,
and $\varepsilon$ is a random error independent of $X$.
Assume that $f,K,\beta_{1},\ldots,\beta_{K}$ are all
unknown, and we observe a random sample $(X_{i},Y_{i}),1\le i\le n$,
which are i.i.d. This is similar to the the setting of Ferr\'e and Yao
(\citeyear{r13}, \citeyear{r14}). Following \citet{r22}, we call
$\beta_1,\ldots, \beta_K$
the EDR directions, and $\operatorname{span}(\beta_1,\ldots,\beta_K)$
the EDR space. Without fixing $f$, the EDR directions are
not identifiable; however, the EDR space is identifiable.
The focus of this paper is the estimation of the dimension, $K$, of the
EDR space.

We assume that the $X_i$'s are observed at each $t\in\CI$.
For functional data, this is an idealized assumption as no functions
on a continuum can be fully observed. However, it is a reasonable
approximation for densely observed smooth data, for which the Tecator
data discussed
in Section \ref{sec:data} is a good example.
In that situation, for most theoretical and practical purposes, one can
fit continuous curves to
the discrete-time data and then treat the fitted curves as the true
functional data;
see, for example, \citet{r17}, \citet{r4} and
\citet{r36}.
The case of sparsely observed functional data requires more attention
and will not be
studied in this paper. It may also be of interest to study the case
where $X$ contains
measurement error; see (a) of Section \ref{s:discussions}.

\subsection{Principal components}\label{s:pca}

First, we focus on the generic process $X$. Denote the mean functions
$\mu$ of $X$ by $\mu(t)=\E\{X(t)\}, t\in\CI$. The covariance
operator of $X$
is the linear operator $\Gamma_{X}:=\E((X-\mu)\otimes(X-\mu))$, where,
for any $h\in\CH$, $h\otimes h$ is the linear operator that maps
any $g\in\CH$ to ${\langle}h,g{\rangle}h$.
It can be seen that $\Gamma_{X}$ is a well-defined compact operator so long
as $\E(\|X\|^4)<\infty$, which we assume throughout the paper; see
\citet{r11a} for the mathematical details in constructing
$\mu$ and $\Gamma_X$. Then there exist nonnegative real
numbers $\omega_1\ge\omega_2\ge\cdots,$ where $\sum_j\omega_j <
\infty$,
and orthonormal functions $\psi_1,\psi_2,\ldots$ in $\CH$ such that
$\Gamma_{X}\psi_j = \omega_j\psi_j$ for all $j$;
namely, the $\omega_j$'s are the eigenvalues and $\psi_j$'s the corresponding
eigenfunctions of $\Gamma_{X}$. The $\psi_j$'s are commonly referred
to as the principal
components of $X$.
%
%
%
It follows that
%
%
\begin{equation}\label{eq:coveigen}
\Gamma_{X}=\sum_{j} \omega_j \psi_j\otimes\psi_j
\end{equation}
and
%
%
\begin{equation}\label{eq:fpca}
X=\mu+\sum_{j} \xi_j \psi_j
= \mu+\sum_{j} \sqrt{\omega_j}\eta_j \psi_j,
\end{equation}
where the $\xi_j$'s are zero-mean, uncorrelated random variables with
$\var(\xi_j)=\omega_j$, and the $\eta_j$'s are standardized $\xi_j$'s.
Call $\eta_j$ the standardized $j$th principal component score of $X$.
The representations in (\ref{eq:coveigen}) and (\ref{eq:fpca})
are commonly referred to as the principal component decomposition and
the Karhunen--Lo\`eve expansion,
respectively; see \citet{r2} and \citet{r11a}
for details.

In view of (\ref{eq:fullmodel}) and (\ref{eq:fpca}),
any component of $\beta_k$ that is in the orthogonal complement of the span
of the $\psi_j$ is not estimable. As explained above, this paper does
not address
the estimation of the $\beta_k$. Thus, assume without generality that the
$\beta_{k}$'s are spanned by the $\psi_{j}$'s and write
%
%
\begin{equation}\label{e:betarep}
\beta_k=\sum_{j} {b_{kj}\over\sqrt{\omega_j}} \psi_j.
\end{equation}
By (\ref{eq:fpca}) and (\ref{e:betarep}),
${\langle}\beta_k,X{\rangle}= {\langle}\beta_k,\mu{\rangle}+
\sum_j b_{kj}\eta_j$, and
(\ref{eq:fullmodel})
can be re-expressed as
\[
Y=f \biggl(\sum_j b_{1j}\eta_j,\ldots, \sum_jb_{Kj}\eta
_j,\varepsilon\biggr),
\]
where, for simplicity, the constants ${\langle}\beta_1,\mu{\rangle
},\ldots,{\langle} \beta_K,\mu{\rangle}$ are absorbed by $f$. For the
i.i.d. sample $(X_1,Y_1), \ldots, (X_n,Y_n)$, let $\eta_{ij}$ be the
standardized $j$th principal component score of $X_i$, and write
%
%
\begin{equation}\label{e:fullmodelpc}
Y_i=f \biggl(\sum_j b_{1j}\eta_{ij},\ldots, \sum_jb_{Kj}\eta
_{ij},\varepsilon_i \biggr).
\end{equation}

\subsection{Elliptically contoured distributions}

As mentioned in Section \ref{sec1}, the relevance of elliptical
symmetry is
evident in the inference of (\ref{eq:fullmodel}). We devote this
subsection to a brief introduction of the notion of elliptically
contoured distribution for Hilbert-space-valued variables.

Let $X$ be as defined in Section \ref{s:pca}. By the assumption $\E(\|
X\|^4)<\infty$, the distribution of $X$ is determined by the (marginal)
distributions of the random variables ${\langle}h,X{\rangle}, h\in\CH$.
Say that $X$ has an elliptically contoured distribution if
%
%
\begin{equation}\label{e:sphere}
\E\bigl(e^{i{\langle}h,X-\mu{\rangle}} \bigr) = \phi({\langle
}h,\Sigma h{\rangle}),\qquad
h\in\CH,
\end{equation}
for some some function $\phi$ on $\bbR$ and self-adjoint, nonnegative
operator $\Sigma$.
Recall that
$X$ is said to be a Gaussian process if ${\langle}h,X{\rangle}$ is normally
distributed for any $h\in\CH$,
and so (\ref{e:sphere}) holds with $\phi(t) = \exp(-t/2)$ and
$\Sigma=\Gamma_{X}$.
However, (\ref{e:sphere}) in general describes a much larger class of
distributions.

The mathematics necessary to characterize elliptically contoured
distributions was
worked out in Schoenberg (\citeyear{sch38}); see Cambanis, Huang and
Simons (\citeyear{chs81})
and Li (\citeyear{r24}). It follows that definition (\ref
{e:sphere}) implies
that $\Sigma$ is a constant
multiple of $\Gamma_{X}$ and $\phi(t^2)$ is a characteristic function.
More explicitly, (\ref{e:sphere})~leads to the characterization
%
%
\begin{equation}\label{e:sph}
X-\mu\eid\Theta\check X,
\end{equation}
where $\Theta$ and $\check X$ are independent, $\Theta$ is a
nonnegative random
variable with $\E(\Theta^2)=1$ and $\check X$ has the same covariance
operator as $X$;
if $X\in\bbR^p$ and $\operatorname{rank}(\Gamma_{X})=k\ge1$ then $\check
X\eid
\Theta A_{p\times k} U_{k\times1}$ where $AA\tran=\Gamma_{X}$ and
$U$ is uniformly distributed on the $k$-dimensional sphere with radius
$\sqrt{k}$;
if $\operatorname{rank}(\Gamma_{X})=\infty$ then
$\check X$ is necessarily a zero-mean Gaussian process.
Recall that $U_{k\times1}$ is asymptotically Gaussian [cf. \citet
{r38}] and
so the infinite-dimensional representation can be viewed as the limit
of the
finite-dimensional one.

\subsection{Functional inverse regression}

To introduce functional inverse regression, we first
state some conditions:
\begin{enumerate}[(C1)]
\item[(C1)]
$\E(\|X\|^4)<\infty$.
\item[(C2)]
For any function $b\in\CH$, there exist some constants
$c_0,\ldots,c_K$ such that
\[
\E({\langle}b,X{\rangle}| {\langle}\beta_1,X {\rangle}, \ldots,
{\langle}\beta_K, X{\rangle})=c_0+c_1 {\langle}\beta_1,X{\rangle
}+\cdots
+c_K{\langle}\beta_K, X{\rangle}.
\]
\item[(C3)]
$X$ has an elliptically contoured distribution;
namely, (\ref{e:sph}) holds.
\end{enumerate}
Conditions (C1)--(C3) are standard conditions in the
inverse regression literature; see, for instance, Ferr\'e and Yao
(\citeyear{r13}, \citeyear{r14}).
As mentioned earlier, condition (C1) guarantees the
principal decomposition;
moreover, it also ensures the convergence rate of
$n^{-1/2}$ in the estimation of the eigenvalues and eigenspaces of
$\Gamma_{X}$
based on an i.i.d. sample $X_1,\ldots,X_n$; see \citet{r10}.
Condition (C2) is a direct extension of
(3.1) in \citet{r22} which addresses multivariate data.
If $X$ is a Gaussian process, then projections of $X$ are jointly
normal, from
which (C2) follows easily.
Condition (C3) describes a broader class of processes
satisfying (C2) than the Gaussian process; for convenience
(C3) is often assumed in lieu of~(C2).

Call the collection $\{\E(X(t)|Y),t\in\CI\}$ of random variables the
inverse regression process and
denote its covariance operator by $\Gamma_{X|Y}$. We will use the notation
$\operatorname{Im}(T)$, for any operator $T$, to denote the range of $T$.
The following result, first appeared in \citet{r13},
is a straightforward extension of Theorem 3.1 of \citet{r22}.
%
\begin{theorem}\label{thm:fsir}
Under \textup{(C1)} and \textup{(C2)}, $\operatorname{Im}(\Gamma_{X|Y})
\subset\operatorname{span}(\Gamma_{X} \beta_1,\ldots,\Gamma_{X}
\beta_K)$.
\end{theorem}

Theorem \ref{thm:fsir} implies that $\operatorname{span}(\Gamma_{X}
\beta
_1,\ldots,
\Gamma_{X} \beta_K)$ contains all of the eigenfunctions that
correspond to the
nonzero eigenvalues of $\Gamma_{X|Y}$.
Consequently, if $\Gamma_{X|Y}$ has $K$ nonzero eigenvalues, then the
space spanned by
the eigenfunctions is precisely $\operatorname{span}(\Gamma_{X} \beta
_1,\ldots
, \Gamma_{X} \beta_K)$.
In that case, one can in principle estimate $\operatorname{span}(\beta
_1,\ldots
,\beta_K)$ through estimating
both $\Gamma_{X}$ and $\Gamma_{X|Y}$. This forms the basis for the
estimation of the
EDR space [cf. \citet{r22} and Ferr\'e and
Yao (\citeyear{r13}, \citeyear{r14})].

While $\Gamma_{X|Y}$ is finite-dimensional under (C1) and (C2),
if $\CH$ is infinite-dimensional then its definition still involves
infinite-dimensional random functions.
In order to implement any inference procedure, we consider a
finite-dimensional adaptation
using principal components.

Let $m$ be any positive integer, where $m\le n-1$ and, if $X$ is
$p$-dimensional,
$m\le p$. Define
$\bdb_{k,(m)}=(b_{k1},\ldots, b_{km})\tran$, $\bdeta_{i,(m)}=(\eta
_{i1},\ldots,\eta_{im})\tran$
and $\varsigma_{ik}=\sum_{j>m} \eta_{ij} b_{kj}$. Then (\ref
{e:fullmodelpc}) can be expressed
as
%
%
\begin{equation}\label{eq:fullmodel1}
Y_i=f\bigl(\bdb_{1,(m)}\tran\bdeta_{i,(m)}+\varsigma_{i1}, \ldots,
\bdb_{K,(m)}\tran\bdeta_{i,(m)}+\varsigma_{iK}, \varepsilon_i\bigr).
\end{equation}
If one regards the $\bdeta_{i,(m)}$ as predictors and combine
the $\varsigma_{ik}$ with $\varepsilon_i$ to form the error, then
(\ref{eq:fullmodel1})
bears considerable similarity with the multivariate model of \citet{r22}.
One fundamental difference is that although the $\varsigma_{ik}$
are uncorrelated with~$\bdeta_{i,(m)}$, they might not be independent of
$\bdeta_{i,(m)}$, unless $X$ is Gaussian.
Another major difference is that we do not directly observe $\bdeta_{i,(m)}$
so that this model might be viewed as a variation of the
errors-in-variables model in
\citet{r7}.
Our estimator for $K$ will be motivated by the finite-dimensional model
(\ref{eq:fullmodel1}). The details of the procedure, including the
role of $m$, will be explained in
Section \ref{sec:tests}. To pave the way for that, we briefly discuss
the inference of the $\bdb_{k,(m)}$ below.

We first need to estimate $\bdeta_{i,(m)}$.
Let
\[
\bar X = n^{-1}\sum_{i=1}^n X_i \quad\mbox{and}\quad\widehat\Gamma_{X}
= n^{-1}\sum_{i=1}^n (X_i-\bar X)\otimes
(X_i-\bar X)
\]
be the sample mean function and the sample covariance operator, respectively.
Let $\widehat\omega_j$ and $\widehat\psi_j$ be the $j$th sample
eigenvalue
and eigenfunction of
$\widehat\Gamma_{X}$. By \citet{r10}, $\widehat
\omega
_j$ and $\widehat\psi_j$ are root-$n$
consistent under (C1). The standardized $j$th principal component
scores of $X_i$
are then estimated by $\widehat\eta_{ij}=\widehat\omega_j^{-1/2}
{\langle}\widehat\psi
_j, X_i-\bar X{\rangle}$; let
$\widehat\bdeta_{i,(m)}=(\widehat\eta_{i1},\ldots,\widehat\eta
_{im})\tran$.

Based on the ``data'' $(\widehat\bdeta_{i,(m)},Y_i), 1\le i\le n$, the
usual sliced inverse
regression (SIR) algorithm can be carried out as follows. Partition the
range of $Y$ into
disjoint intervals, $S_h$, $h=1,\ldots, H$, where $p_h:=\bbP(Y\in
S_h)>0$ for all $h$.
Define
%
%
\begin{equation}\label{eq:within_slice_mean_vect}
\vartheta_{j,h}=\E(\eta_j|Y\in S_h),\qquad \bdvartheta_{h,(m)}=\E
\bigl(\bdeta_{(m)}|Y\in
S_h\bigr)=(\vartheta_{1,h},\ldots,\vartheta_{m,h})\tran\hspace*{-26pt}
\end{equation}
and
%
%
\begin{equation}\label{eq:within_slice_mean_curve}
\mu_h = \E(X|Y\in S_h)=\mu+\sum_{j}\omega_j^{1/2} \vartheta_{j,h}
\psi_j.
\end{equation}


Let $V_{(m)}=\sum_h p_h \bdvartheta_{h,(m)} \bdvartheta_{h,(m)}\tran
$ be the
between-slice covariance matrix. In the finite-dimensional model (\ref
{eq:fullmodel1}) with
$(\varsigma_{i1},\ldots,\varsigma_{iK}, \varepsilon_i)$ playing the
role of error,
$V_{(m)}$~is the sliced-inverse-regression covariance
matrix. The eigenvectors of $V_{(m)}$ corresponding to the nonzero
eigenvalues are contained\vspace*{2pt} in
$\operatorname{span}(\bdb_{1,(m)}, \ldots,\bdb_{K,(m)})$. The matrix
$V_{(m)}$ is estimated by the corresponding sample version
$\widehat V_{(m)}=\sum_{h=1}^{H}\widehat p_h \widehat\bdvartheta
_{h,(m)}(\widehat
\bdvartheta_{h,(m)})\tran$, where
%
%
\begin{eqnarray}\label{e:phnh}
\widehat p_h&=&n_h/n,\qquad n_h=\sum_i I(Y_i\in S_h) \quad\mbox
{and}\nonumber\\[-8pt]\\[-8pt]
\widehat\bdvartheta_{h,(m)}&=&{1\over n_h} \sum_{i=1}^{n} \widehat
\bdeta
_{i,(m)} I(Y_i \in S_h).\nonumber
\end{eqnarray}
Letting $\widehat\bdb_{1,(m)}, \ldots, \widehat\bdb_{K,(m)}$ be
the first
$K$ eigenvectors of $\widehat V_{(m)}$,
the estimators of $\beta_k$'s are given by
\[
\widehat\beta_k(t)=\sum_{j=1}^m \widehat\omega_j^{-1/2} \widehat
b_{kj} \widehat\psi
_j(t),\qquad k=1,\ldots, K.
\]
In order for $\operatorname{span}(\widehat\beta_1,\ldots,\widehat\beta_K)$ to consistently
estimate the EDR space,
it is necessary that $\operatorname{span}(\mu_1-\mu,\ldots,\mu_h-\mu)$ have the
same dimension as the EDR space,
and that $m$ tends to $\infty$ with $n$ in some manner.
However, first and foremost, we must know $K$ beforehand,
which makes the determination of $K$ a fundamental issue.

The matrix $\widehat V_{(m)}$ will be our basis for deciding $K$.
Here, we define some notation related to $V_{(m)}$ and $\widehat V_{(m)}$
for future use.
For any $m\times1$ vector
$\bdu$, let $\CJ_\bdu=I-\bdu\bdu\tran$; let $\bdg=(g_1,\ldots
,g_H) = (p_1^{1/2},\ldots,p_H^{1/2})$,
and $\widehat\bdg=(\widehat g_1,\ldots, \widehat g_H) = (\widehat
p_1^{1/2},\ldots,\widehat
p_H^{1/2})$.
Define
\begin{eqnarray*}
M&=&\bigl[\bdvartheta_{1,(m)},\ldots,\bdvartheta_{H,(m)}\bigr]_{m\times
H},\qquad G=\diag\{g_1,\ldots, g_H\},\\
F&=&G\CJ_\bdg,\qquad B_{(m)}=M F, \\
\widehat M &=& \bigl[\widehat\bdvartheta_{1,(m)}, \ldots,
\widehat\bdvartheta_{H,(m)}\bigr]_{m\times H},\qquad \widehat G =\diag\{
\widehat
g_1,\ldots, \widehat g_H\},\\
\widehat F&=&\widehat G \CJ_{\widehat\bdg}, \qquad
\widehat B_{(m)}=\widehat M \widehat F,
\end{eqnarray*}
where $\bdvartheta_{h,(m)}$ and $\widehat\bdvartheta_{h,(m)}$
are defined in (\ref{eq:within_slice_mean_vect}) and (\ref{e:phnh}),
respectively. Thus, the
inverse-regression covariance matrices $V_{(m)}$ and $\widehat V_{(m)}$ can
be rewritten as
%
%
\begin{equation}\label{e:VB}
V_{(m)}=B_{(m)} B_{(m)}\tran,\qquad \widehat V_{(m)}=\widehat B_{(m)}
\widehat
B_{(m)}\tran.
\end{equation}

%
%
%
%
%
%
%

%
\section{Deciding the dimension of EDR space} \label{sec:tests}

As explained in previous sections, we are particularly interested in
functional data or high-dimensional multivariate data. Existing
methods for deciding the dimensionality of EDR space in the multivariate
setting [\citet{r22}, \citet{r30}] are not directly
applicable to the types of
data that are focused on in this paper. Ferr\'e and Yao (\citeyear
{r13}, \citeyear{r14})
used a
graphical approach to determine the number of EDR directions for
functional data
but a formal statistical procedure has been lacking.

Our approach is generically described as follows.
To decide the dimension of the EDR space, as in \citet{r22}, we will
conduct sequential testing of $H_0\dvtx K\le K_0$ versus $H_a\dvtx K >
K_0$ for
$K_0=0,1,2,\ldots;$
we will stop at the first instance $K_0=\widehat K$ when the test fails to
reject $H_0$
and declare $\widehat K$ as the true dimension.
Below, we consider two types of tests in the sequential testing procedure
motivated by (\ref{eq:fullmodel1}).
In Section \ref{s:chisq}, we assume that $m$ is fixed, while in Section
\ref{sec:Neyman} we
consider $m$ in a wide range.

\subsection{Chi-squared tests based on a fixed $m$}\label{s:chisq}

Fix an $m$ and focus on the between-slice inverse covariance matrix
$V_{(m)}$, which
has dimension $m\times m$; recall that it only makes sense to consider
$m$ such that
$m\le n-1$ and, if $X$ is a $p$-dimensional vector, $m\le p$.
Define
\[
K_{(m)}=\operatorname{rank}\bigl(V_{(m)}\bigr).
\]
Clearly, $K_{(m)}\le K$ for all $m$.
It is desirable to pick an $m$ such that $K_{(m)}=K$. Note that this
condition means that
the projections of all of the EDR directions onto the space spanned by
the first $m$
principal components are linearly independent, which is very different
from saying that
all of the EDR directions are completely in the span of the first $m$
principle components; see the examples in Section \ref{sec:simulation}.
However, picking an $m$ to guarantee $K_{(m)}=K$ before analyzing the data
is clearly not always possible. A practical approach is to simply pick
an $m$ such that the first
$m$ principal components explain a large proportion, say, $95\%$, of
the total
variation in the $X_i$'s. Such an approach will work for most
real-world applications.
Still, keeping $m$ fixed has its limitations. We will address them in
more detail in future sections.

In the following, let $\lambda_j(M)$ denotes the $j$th largest
eigenvalue of a\break nonnegative-definite
square matrix $M$. Under $H_0\dvtx K\le K_0$, we have\break
$\lambda_{K_0+1}(V_{(m)})=\cdots=\lambda_m(V_{(m)})=0$.
Consider the test statistic
%
%
\begin{equation}\label{eq:aveigval}
\CT_{K_0,(m)}=n \sum_{j=K_0+1}^m \lambda_j\bigl(\widehat V_{(m)}\bigr).
\end{equation}
Since $\widehat V_{(m)}$ estimates $V_{(m)}$, large values of $\CT_{K_0,(m)}$
will support the rejection of~$H_0$.
The following theorem provides the asymptotic distribution of
$\CT_{K_0,(m)}$ under~$H_0$.
For the convenience of the proofs, we will assume below that
the positive eigenvalues of $\Gamma_{X}$ are all distinct.

The following addresses the case where $X$ is a Gaussian process.
\begin{theorem}\label{thm:fsir_asymp}
Suppose that \textup{(C1)} holds and $X$ is a Gaussian process.
Assume that $K\le K_0$, and let $H>K_0+1$ and $m\ge K_0+1$.
Denote by $\CX$ a random variable having a $\chi^2$ distribution with
$(m-K_0)(H-K_0-1)$ degrees of freedom.
\begin{longlist}
\item
If $K_{(m)}=K_0$, then
%
%
\begin{equation}\label{eq:eigasymp}
\CT_{K_0,(m)}\cid\CX\qquad\mbox{as $n\to\infty$}.
\end{equation}
\item
If $K_{(m)}<K_0$, then $\CT_{K_0,(m)}$ is asymptotically stochastically
bounded by $\CX$; namely,
\[
\limsup_{n\to\infty} {\mathbb P}\bigl(\CT_{K_0,(m)} > x\bigr) \le{\mathbb
P}(\CX> x) \qquad\mbox{for all $x$}.
\]
\end{longlist}
\end{theorem}

Theorem \ref{thm:fsir_asymp} suggests a $\chi^2$ test for testing
$H_0\dvtx K\le K_0$ versus $H_a\dvtx K>K_0$,
which is an extension of a test in \citet{r22} for multivariate data.
Ideally, case (i) holds and the $\chi^2$ test has the
correct size asymptotically, as $n\to\infty$.
For a variety of reasons case (ii) may be true, for which the $\chi^2$ test
will be conservative. This point will be illustrated graphically by a
simulation example in
Figure \ref{fig:testmean} in Section~\ref{sec:simulation}.

The proof of Theorem \ref{thm:fsir_asymp} is highly nontrivial, which
goes considerably beyond the scope
of the multivariate counterpart.
A theoretical result that is needed to establish (ii) of
Theorem \ref{thm:fsir_asymp} appears to be new and is stated here.
\begin{pp}\label{pp:eign_value_bound}
Let $Z$ be a $p\times q$ random matrix and we write
$Z=[ Z_1 | Z_2 ]$ where $Z_1$ and $Z_2$
have sizes $p\times r$ and $p\times(q-r)$, respectively, for some $0 <
r < \min(p,q)$.
Assume that $Z_1$ and $Z_2$ are
independent, and $Z_2$ contains i.i.d. $\Normal(0,1)$ entries. Then
$\sum_{j=r+1}^{p} \lambda_j(Z Z\tran)$
is stochastically bounded by $\chi^2$ with $(p-r)(q-r)$ degrees of freedom.
\end{pp}

The case where $Z$ is a matrix of i.i.d. $\Normal(0,1)$ entries can be
viewed as
the special case, $r=0$, in Proposition \ref{pp:eign_value_bound}.
In that case, the bound is the exact distribution
since $\sum_{j=1}^{p} \lambda_j(Z Z\tran)$ equals the sum of squares
of all of the entries
of $Z$ and is therefore distributed as $\chi^2$ with $pq$
degrees of freedom.

Next, we address the scenario where $X$ is elliptically contoured
but not necessarily Gaussian.
Let
%
%
\begin{equation}\label{e:tau}
\tau_h = \E(\Theta^2|Y\in S_h),\qquad h=1,\ldots,H.
\end{equation}
If $K_{(m)}=K_0$, then it can be seen from the proofs in the \hyperref
[app]{Appendix} that
%
%
\begin{equation}\label{eq:eigasymp2}
\CT_{K_0,(m)}\cid\sum_{k=1}^{H-K_0-1} \delta_k
\CX_k \qquad\mbox{as $n\to\infty$},
\end{equation}
where $\CX_k$'s are distributed as i.i.d. $\chi^2$ with $m-K_0$
degrees of freedom,
and $\delta_1,\ldots,\delta_{H-K_0-1}$ are the eigenvalues of
$\Lambda\Xi\Lambda$, with $\Xi=\CJ_{\bdg} \{I-B_{(m)}\tran
(B_{(m)}\times\break B_{(m)}\tran)^-
B_{(m)}\}\CJ_{\bdg}$ and $\Lambda=\diag(\tau_1^{1/2},\ldots,\tau
_H^{1/2})$.
If $X$ is Gaussian, then $\tau_h$'s and $\delta_k$'s are identically
equal to 1.
In general, the limiting null distribution in (\ref{eq:eigasymp2})
depends on the unknown parameters
$\delta_k$. \citet{r9} suggested carrying out this type of test by
simulating the critical regions based on
the estimated values of these parameters. Below, we introduce a
different approach by
adjusting the test statistic so that the limiting distribution is
free of nuisance parameters.

Under $H_0\dvtx K\le K_0$, let $m>K_0$ and $\widehat\CP_2$ be the matrix
whose columns
are the eigenvectors that correspond to the $m-K_0$
smallest eigenvalues of $\widehat V_{(m)}$. The definition (\ref
{e:tau}) suggests
(see proof of Theorem \ref{thm:fsir_asymp1} in the \hyperref[app]{Appendix})
that $\tau_h$ can be estimated by
%
%
\begin{eqnarray}\label{eq:cond_var_est}
\widehat\tau_h&=&{1\over(m-K_0) n_h} \tr\Biggl\{\widehat\CP_2
\widehat{\CP}
{\,}_{2}^{\mathrm{T}}\sum_{i=1}^n
\bigl(\widehat\bdeta_{i,(m)}-\widehat\bdvartheta_{h,(m)}\bigr)\nonumber\\[-8pt]\\[-8pt]
&&\hspace*{119.8pt}{}\times
\bigl(\widehat\bdeta_{i,(m)}-\widehat\bdvartheta_{h,(m)}\bigr)\tran I(Y_i\in
S_h) \Biggr\}.\nonumber
\end{eqnarray}
%
%
%
%
%
Put $\Lambda=\diag(\tau_1^{1/2},\ldots, \tau_H^{1/2})$,
$\widehat\Lambda=\diag(\widehat\tau_1^{1/2},\ldots, \widehat\tau
_H^{1/2})$, and
define
\begin{eqnarray*}
W_{(m)}&=&B_{(m)} \Lambda(\Lambda\CJ_{\bdg} \Lambda)^-,\qquad
\Sigma_{(m)}= W_{(m)} W_{(m)}\tran, \\
\widehat W_{(m)}&=&\widehat B_{(m)} \widehat\Lambda(\widehat
\Lambda\CJ_{\widehat\bdg}
\widehat\Lambda)^{-},\qquad
\widehat\Sigma_{(m)}=\widehat W_{(m)}
\widehat W_{(m)}\tran,
\end{eqnarray*}
where $A^-$ denotes the Moore--Penrose generalized inverse of the matrix
$A$. By Lemma \ref{lm:columnspace_W} below,
$\Sigma_{(m)}$ has the same null space as $V_{(m)}$. Thus, under
$H_0\dvtx
K\le K_0$,
we have $\lambda_{K_0+1}(\Sigma_{(m)})=\cdots=\lambda_m(\Sigma
_{(m)})=0$. We therefore
propose the test statistic
\[
\CT_{K_0,(m)}^{\ast}=n\sum_{j=K_0+1}^m \lambda_j\bigl(\widehat\Sigma_{(m)}\bigr),
\]
where, again,\vspace*{2pt} large values of $\CT_{K_0,(m)}^\ast$ support the
rejection of $H_0$.
The following result extends (i) of Theorem \ref{thm:fsir_asymp} from
the case where $X(t)$ is
Gaussian to a general elliptically contoured process. While we
conjecture that (ii) of Theorem \ref{thm:fsir_asymp}
can be similarly extended, we have not been able to prove it.
\begin{theorem}\label{thm:fsir_asymp1}
Suppose that \textup{(C1)} and \textup{(C3)} hold.
Assume that the true dimension $K\le K_0$ and let $H>K_0+1$ and $m\ge K_0+1$.
If $K_{(m)}=K_0$ then $\CT_{K_0,(m)}^{\ast}\cid\chi
^2_{(m-K_0)(H-K_0-1)}$ as $n\to\infty$.
\end{theorem}

The test of $H_{0}\dvtx K\le K_{0}$ based on
$\CT_{K_0,(m)}$ and $\CT_{K_0,(m)}^{\ast}$ and the
asymptotic null distribution, $\chi^2_{(m-K_0)(H-K_0-1)}$, will be
referred to
as the $\chi^2$ test and the adjusted $\chi^2$ test, respectively.

%
\subsection{Adaptive Neyman tests}
\label{sec:Neyman}
So far we considered tests based on a fixed $m$. In most situations, in
practice, choosing the smallest
$m$ for which the first $m$ sample principal components explain most of the
variations should work fairly well for determining $K$. However, for
functional or high-dimensional
multivariate data, one cannot theoretically rule out the possibility
that the EDR directions can only be
detected by examining the information contained in higher-order
principal components.

A careful inspection reveals two different issues here. The first is
the question that if we have
an unusual model in which some EDR directions depend only on high-order
principal components
that the data have little power to discern, can any approach be
effective in detecting the
presence of those directions? The answer is ``not likely'' since,
intuitively, we can detect the
presence of those EDR directions no better than we can the principle
components that
comprise the directions. This is due more to the nature of high or
infinite-dimensional data
than the limitation of any methodology. However, keep in mind that
principal components
are not ordinary covariates, but are mathematical constructs which not
only depend on
the covariance function of $X$ but also the choice of inner product of
$\CH$.
Thus, one can argue that having an EDR direction that is
orthogonal to a large number of low-order principal components of the
predictor is itself a
rather artificial scenario and is not likely to be the case in practice.


Let us now turn to the second issue. Assume that all of the EDR
directions do contain low-order principal components which can be
estimated well from data.
For example, suppose each EDR direction is not in the orthogonal
complement of the space
spanned by the first three principal components and so
the procedures described in Section \ref{s:chisq} will in principle
work if we
let $m=3$. However, since
that knowledge is not available when we conduct data analysis, to be
sure perhaps we
might consider picking a much larger truncation point, say, $m=30$.
The problem with this approach is that, when the sample size is fixed,
the power of the tests will decrease with $m$.
Intuitively, when $m$ is large the information contained in the
individual components of
$\widehat\bdvartheta_{h,(m)} = (\widehat\vartheta_{1,h},\ldots,
\widehat
\vartheta_{m,h})\tran$
becomes diluted. We will illustrate this point numerically in Section
\ref{sec:simulation1}.
This is strikingly similar to the situation of testing whether the mean
of a
high-dimensional normal random vector is nonzero described at the
beginning of Section 2
of \citet{r12}, where the power of the Neyman test
(likelihood-ratio test) was shown to
converge to the size of the test as the number of dimension increases.
Essentially, the problem
that they describe is caused
by the fact that the Neyman test has a rejection region that is
symmetric in all components of the vector,
which is designed to treat all possible alternatives uniformly.
\citet{r12}
argued that the alternatives that are of the most importance in
practice are usually those in
which the leading components of the Gaussian mean vector are nonzero,
and they modified the
Neyman test accordingly such that the test will have satisfactory
powers for those alternatives.

We now introduce a test inspired by \citet{r12} that avoids
having to
pick a specific $m$. To test $H_0\dvtx K\le K_0$ against $H_a\dvtx K>K_0$,
we compute $\CT_{K_0,(m)}$ for all of $m=K_0+1,\ldots, N$, for some
``large'' $N$; we
then take the maximum of the standardized versions of these test statistics,
and the null hypothesis will be rejected for large values of the maximum.
To facilitate this approach, we present the following result that is a
deeper version of
Theorem \ref{thm:fsir_asymp} and shows that the test statistics $\CT
_{K_0,(m)}$ has a
``partial sum'' structure in $m$ as the sample size tends to $\infty$.
\begin{theorem}\label{thm:increasing_m}
Suppose that \textup{(C1)} holds and $X$ is a Gaussian process.
Assume that $K\le K_0$ and let $H>K_0+1$.
Let $\chi_i^2, i\ge1$, be i.i.d. $\chi^2$ random variables with
$H-K_0-1$ degrees of freedom and
define $\CX_{(m)} = \sum_{i=1}^{m-K_0} \chi_i, m\ge K_0+1$.
Then, for all positive integers $N>K_0$, the collection of test
statistics $\CT_{K_0,(m)}, m = K_0+1,\ldots, N$, are jointly
stochastically bounded by
$\CX_{(m)}, m = K_0+1,\ldots, N$, as the sample size $n$ tends to
$\infty$.
%
%
%
\end{theorem}

In view of Theorem \ref{thm:increasing_m}, we propose the following.
To test $H_0\dvtx K\le K_0$ versus $H_a\dvtx K>K_0$, define
\[
\CU_{K_0,N}:= \max_{K_0+1\le m \le N}
{ \CT_{K_0,(m)}-(m-K_0)(H-K_0-1) \over
\sqrt{2(m-K_0)(H-K_0-1)}},
\]
and we reject $H_0$ at level $\alpha$ if $\CU_{K_0,N}>u_{\alpha}$
where $u_{\alpha}$ is the
$1-\alpha$ quantile of
\[
\CB_{K_0,N}:=\max_{K_0+1\le m \le N}
{ \CX_{(m)}-(m-K_0)(H-K_0-1) \over
\sqrt{2(m-K_0)(H-K_0-1)}}.
\]
The resulting test resembles asymptotically the aforementioned test in
\citet{r12}
which was referred to as an adaptive Neyman test. For convenience,
we will also refer to our test as adaptive Neyman test, although the
contexts of the
two problems are completely unrelated.

Suppose that $H_0$ holds and $m_{K_0}$ is the smallest $m$ such that
$K_{(m)}=K_0$.
Then, by Theorem \ref{thm:fsir_asymp}, $\CU_{K_0,N}-\CU_{K_0,m_{K_0}-1}
\cid\CB_{K_0,N}-\CB_{K_0,m_{K_0}-1}$.
Thus, $\CB_{K_0,N}$ is, intuitively, a tight asymptotic stochastic
bound for $\CU_{K_0,N}$.

Simulation results show that the maximum in the definition of $\CU
_{K_0,N}$ is, with high probability,
attained at relatively small $m$'s. Thus, the test is quite robust with
respect to the choice of $N$.
In practice, one can pick $N$ so that there is virtually just noise
beyond the $N$th sample principal component.
Numerically, the performance of the adaptive Neyman test matches those
of the $\chi^2$ tests
in which $m$ is chosen correctly, but does not have the weakness of
possibly under-selecting $m$.

\subsection{Discussion}\label{s:discussions}

\begin{enumerate}[(a)]
\item[(a)]
Our procedures apply to both finite-dimensional and
infinite-dimensional data, and,
in particular, are useful for treating high-dimensional multivariate data.
In that case, Li's $\chi^2$ test suffers from the problem of diminishing
power as does the test developed in \citet{r30}; see, for example,
Table 4 in \citet{r30}.
Our procedures can potentially be a viable solution in overcoming
the power loss problem in that situation. The inclusion of measurement
error in $X$
provides additional flexibility in modeling multivariate data. Note
that the formulation of
Theorem \ref{thm:fsir} can be extended to accommodate measurement error:
if $X=X_1+X_2$ where $X_1$ is the true covariate with mean $\mu$ and covariance
matrix $\Gamma_{X_1}$ and $X_2$ is independent measurement error with
mean zero, then
$\E(X|Y)=\E(X_1|Y)$ and so $\operatorname{Im}(\Gamma_{X|Y})
\subset\operatorname{span}(\Gamma_{X_1} \beta_1,\ldots,\Gamma_{X_1}
\beta_K)$.
Thus, one might speculate that our procedures continue to work in that
case, and this is borne out by simulations presented in Section \ref
{s:multivariate}. Detailed theoretical investigation of this is a topic
of future work, but preliminary indications are that the extension of
Theorem \ref{thm:fsir} is valid at least under the additional
assumption that the components of $X_2$ are i.i.d. with finite variance.
%
\item[(b)]
Choice of slices: in the SIR literature, the prevailing view is that
the choice of slices is of
secondary importance. In our simulation studies, we used contiguous
slices containing roughly the same
number of $Y_i$'s, where the number of $Y_i$'s per slice that we
experimented with ranged from
25 to 65. Within this range, we found that the number of data per slice
indeed had a negligible effect on
the estimation of $K$.
\item[(c)]
Choice of $\alpha$: if $\alpha$ is fixed and $m$ and $N$ are chosen
sensibly in the
$\chi^2$ tests and the adaptive Neyman test, respectively, then the
asymptotic results show that the probability of correct identification
of $K$ tends to $1-\alpha$ as $n$ tends to~$\infty$.
In real-data applications, the optimal choice of $\alpha$ depends on a
number of factors
including the sample size and the true model. In our simulation
studies, presented in Section \ref{sec:simulation},
$\alpha=0.05$ worked well for all of our settings.
\item[(d)]
Limitations of SIR: the failure of SIR in estimating the EDR space in
situations where $Y$
depends on $X$ in a symmetric manner is well documented. While exact
symmetry is not a
highly probable scenario in practice,
it does represent an imperfection of SIR which has been addressed by a
number of other methods including
SAVE in \citet{r8}, MAVE in \citet{r32} and
integral transform
methods in Zhu and Zeng (\citeyear{r34}, \citeyear{r35}). The
estimation of $K$ based on
those approaches will
be a topic of future research.
\end{enumerate}

%
\section{Simulation studies}
\label{sec:simulation}
\subsection{Simulation 1: Sizes and power of the
tests}\label{sec:simulation1} In this study, we consider functional
data generated
from the process
\[
X(t)=\sum_{k=1}^\infty\omega_{2k-1}^{1/2} \eta_{2k-1}
\sqrt{2}\cos(2k\pi t)+\sum_{k=1}^\infty\omega_{2k}^{1/2} \eta_{2k}
\sqrt{2}\sin(2k\pi t),
\]
where $\omega_k=20 (k+1.5)^{-3}$. Thus, the principal
components are the sine and cosine curves in the sum. We will consider
the cases
where the $\eta_k$'s follow $\Normal(0,1)$ and centered
multivariate $t$ distribution with $\nu=5$ degrees of freedom, with
the latter representing the situation where $X$ is an elliptically
contoured process. Note that the centered multivariate $t$
distribution with $\nu$ degrees of freedom can be represented by
\[
\bdt=Z/\sqrt{\tau/(\nu-2)}\sim t_\nu\qquad\mbox{where } Z\sim
N(0,I) \mbox{ and } \tau\sim\chi^2_\nu\mbox{ are independent}.
\]
To simulate $\{\eta_1,\eta_2,\ldots\}$ in that case, we first
simulate $z_1,z_2\cdots\sim\mathrm{i.i.d.} \Normal(0,\break1)$, $\tau\sim
\chi^2_\nu$, and then put $\eta_k=z_k/\sqrt{\tau/(\nu-2)}$. By this
construction, any finite collection of the $\eta_k$'s follows a
multivariate $t$ distribution, where the $\eta_k$'s are mutually
uncorrelated but not independent.


Let the EDR space be generated by the functions
%
%
\begin{eqnarray}\label{e:dir}\qquad\quad
\beta_1(t)&=&0.9\sqrt{2}\cos(2\pi t)+1.2\sqrt{2}\cos(4\pi
t)\nonumber\\
&&{}-0.5\sqrt{2}\cos(8\pi
t)+\sum_{k>4} {\sqrt{2}\over(2k-1)^3} \cos(2k\pi t), \nonumber\\
\beta_2(t)&=&-0.4 \sqrt{2}\sin(2\pi t)+ 1.5 \sqrt{2}\sin(4\pi
t)-0.3 \sqrt{2}\sin(6\pi
t) \nonumber\\[-8pt]\\[-8pt]
&&{}+0.2\sqrt{2}\sin(8\pi t)+\sum_{k>4} {(-1)^k \sqrt{2} \over
(2k)^3} \sin(2k\pi
t),\nonumber\\
\beta_3(t)&=&\sqrt{2}\cos(2\pi t)+\sqrt{2}\sin(4\pi t)+0.5\sqrt
{2}\cos( 6\pi t)+0.5\sqrt{2} \sin
(8\pi t) \nonumber\\
&&{}+\sum_{k\ge3} {\sqrt{2}\over(4k-3)^3} \cos\{2(2k-1)\pi t\}
+\sum_{k\ge3}{\sqrt{2}\over(4k)^3}\sin(4k\pi t). \nonumber
\end{eqnarray}
Consider the models
\begin{eqnarray*}
\mbox{Model 1: } Y&=&1+2\sin({\langle}\beta_1, X{\rangle
})+\varepsilon,\\
\mbox{Model 2: } Y&=&{\langle}\beta_1, X{\rangle}\times(2{\langle
}\beta_2,X{\rangle}
+1)+\varepsilon,\\
\mbox{Model 3: } Y&=&5 {\langle}\beta_1, X{\rangle}\times(2{\langle}
\beta_2,X{\rangle}+1)/(1+{\langle}\beta_3,X{\rangle
}^2)+\varepsilon,
\end{eqnarray*}
where $\varepsilon\sim\Normal(0,0.5^2)$. The EDR spaces of the three
models have dimensions $K=1,2$ and $3$, respectively. Also note that
$K_{(m)}=K$ if $m\ge1, 2$ and $3$, respectively, for the three models.

In each of 1000 simulation runs, data were simulated from models 1--3
for the two distributional scenarios that $X$ is distributed as
Gaussian and $t$ with two sample sizes $n=200$ and $500$. To mimic
real applications, we assumed that each curve $X_i$ is observed at
$501$ equally-spaced points. We then registered the curves using 100
Fourier basis functions. A functional principal component analysis
was carried out using the package \texttt{fda} in \texttt{R} contributed
by Jim Ramsay.

To decide the dimension of the EDR space, we compared the two
proposed $\chi^2$ tests and the adaptive Neyman test. For the
$\chi^2$ tests, we let $m=5, 7$ and 30, where the first 5, 7 and 30
principal components of $X$, respectively, account for
$91\%$, $95\%$ and $99.59\%$ of the total variation.
We present the results for $m=30$
as an extreme case to illustrate the point that using a large number of
principal components
will cause the tests to have lower powers. For the adaptive Neyman
test, we took $N=K_0+30$ and simulated the critical values for $\CU_{K_0,N}$
based on the description following Theorem \ref{thm:increasing_m}.
We only report the results based on $H=8$ slices, but the choice was
not crucial. The nominal size of the tests was set to be $\alpha=0.05$.

%
\begin{table}
\caption{Empirical frequencies of rejecting the
hypothesis $H_0\dvtx K\le1$. The results are based on 1000
simulations for each of the three models, two distributions of process
$X(t)$, and two sample sizes. The $\chi^2$ test and the adjusted
$\chi^2$ are applied with fixed $m$ values, and the adaptive Neyman
test is
applied with $N=K_0+30$}
\label{tab:table1}
\begin{tabular*}{\tablewidth}{@{\extracolsep{\fill}}lccccccc@{}}
\hline
& & \multicolumn{2}{c}{\textbf{Model 1}} &
\multicolumn{2}{c}{\textbf{Model 2}} &\multicolumn{2}{c@{}}{\textbf{Model 3}}\\[-4pt]
& & \multicolumn{2}{c}{\hrulefill} &
\multicolumn{2}{c}{\hrulefill} &\multicolumn{2}{c@{}}{\hrulefill}\\
& \textbf{Distribution of} $\bolds\eta$ & \textbf{Normal} & $\bolds t$ & \textbf{Normal}
& $\bolds t$ & \textbf{Normal} & $\bolds t$\\
\hline
$n=200$ & $\chi^2$ test ($m=5$) & 0.040 & 0.046 &0.883 & 0.567 & 0.995
& 0.949\\
& Adj. $\chi^2$ ($m=5$) & 0.070 & 0.082 &0.894 & 0.584 & 0.997 & 0.954
\\
[2pt]
& $\chi^2$ test ($m=7$) &0.045 &0.039 &0.827 & 0.496 &0.982 & 0.910\\
& Adj. $\chi^2$ ($m=7$) &0.071 & 0.083 &0.868 & 0.537 &0.989 & 0.924\\
[2pt]
& $\chi^2$ test ($m=30$) & 0.034 & 0.020 & 0.320 & 0.111 & 0.677 &
0.493\\
& Adj. $\chi^2$ ($m=30$) & 0.105 & 0.072 & 0.484 & 0.299 & 0.798 &
0.668\\
[2pt]
& Adaptive Neyman &0.044 & 0.045 &0.860 & 0.545& 0.993 & 0.936 \\
[4pt]
$n=500$ & $\chi^2$ test ($m=5$) & 0.052 & 0.043 &1.000 & 0.977 & 1.000
& 1.000\\
& Adj. $\chi^2$ ($m=5$) & 0.069 & 0.056 &1.000 & 0.969 & 1.000 &
1.000\\
[2pt]
& $\chi^2$ test ($m=7$) & 0.048 & 0.033 &1.000 & 0.958 & 1.000 &
0.999\\
& Adj. $\chi^2$ ($m=7$) & 0.056 &0.049 &1.000 & 0.949 & 1.000 & 0.999\\
[2pt]
& $\chi^2$ test ($m=30$) & 0.059 & 0.025 & 0.958 & 0.584 & 0.999 &
0.990 \\
& Adj. $\chi^2$ ($m=30$) & 0.085 & 0.054 & 0.972 & 0.666 & 0.999 &
0.991\\
[2pt]
& Adaptive Neyman & 0.052 & 0.040 &1.000 & 0.963 & 1.000 & 0.999 \\
%
%
%
\hline
\end{tabular*}
\end{table}

The simulation results are briefly discussed below.
Table \ref{tab:table1} gives the empirical frequencies of rejecting
$H_0\dvtx K\le1$. Since the dimension of EDR space
under model 1 is equal to 1, the results in the column under model 1
give the empirical sizes of the tests. Models 2 and 3 represent two
cases under the alternative
hypothesis, therefore the results in those columns give the power of
the tests.
As can be seen, when $m$ is $5$ or $7$, the two $\chi^2$ tests have
sizes close to the nominal size
and have high powers. However, for the case $m=30$ and $n=200$, the
tests performed significantly
worse in those metrics. On the other hand the adaptive Neyman test
performs very stably, with powers comparable to those of $\chi^2$
tests with a well chosen $m$.
It is also worth noting that when $X$ has the $t$ distribution,
the $\chi^2$ test performs comparably to, sometimes better than, the
adjusted $\chi^2$ test,
showing that the $\chi^2$ test is quite robust against departure from
normality.



%
\begin{table}
\caption{Empirical frequencies of finding the true
dimension of the EDR space. The results are based on 1000
simulations for each of the three models, two distributions of process
$X(t)$, and two sample sizes. The $\chi^2$ test and the adjusted
$\chi^2$ are applied with fixed $m$ values, and the adaptive Neyman
test is
applied with $N=K_0+30$}
\label{tab:table2}
\begin{tabular*}{\tablewidth}{@{\extracolsep{\fill}}lccccccc@{}}
\hline
& & \multicolumn{2}{c}{\textbf{Model 1}} &
\multicolumn{2}{c}{\textbf{Model 2}} &\multicolumn{2}{c@{}}{\textbf{Model 3}}\\[-4pt]
& & \multicolumn{2}{c}{\hrulefill} &
\multicolumn{2}{c}{\hrulefill} &\multicolumn{2}{c@{}}{\hrulefill}\\
& \textbf{Distribution of} $\bolds\eta$ & \textbf{Normal} & $\bolds t$
& \textbf{Normal} & $\bolds t$ & \textbf{Normal} & $\bolds t$\\
\hline
$n=200$ & $\chi^2$ test ($m=5$) & 0.960 &0.954 & 0.859 & 0.562 & 0.322
& 0.165\\
& Adj. $\chi^2$ ($m=5$) & 0.930 &0.918 & 0.846 & 0.556 & 0.393 &
0.224\\
[2pt]
& $\chi^2$ test ($m=7$) &0.955 & 0.961 & 0.809& 0.488 &0.272 & 0.116\\
& Adj. $\chi^2$ ($m=7$) &0.929 & 0.917 & 0.832& 0.505 & 0.313 & 0.167
\\
[2pt]
& $\chi^2$ test ($m=30$) &0.966 & 0.971 & 0.309 & 0.111 & 0.057 &
0.026 \\
& Adj. $\chi^2$ ($m=30$) &0.895 & 0.924 & 0.462 & 0.277 & 0.119 &
0.079 \\
[2pt]
& Adaptive Neyman &0.956 & 0.955 & 0.843 & 0.542 & 0.337 &0.158 \\
[4pt]
$n=500$ & $\chi^2$ test ($m=5$) & 0.948 &0.957 & 0.955 & 0.958 & 0.842
& 0.629\\
& Adj. $\chi^2$ ($m=5$) & 0.931 &0.944 & 0.948 & 0.910 & 0.849 &
0.648\\
[2pt]
& $\chi^2$ test ($m=7$) & 0.952 &0.967 & 0.959 & 0.948 & 0.739 &
0.489\\
& Adj. $\chi^2$ ($m=7$) & 0.944 &0.951 & 0.949 & 0.898 & 0.754 &
0.551\\
[2pt]
& $\chi^2$ test ($m=30$) & 0.941 & 0.975 & 0.913 & 0.582 & 0.279 &
0.101 \\
& Adj. $\chi^2$ ($m=30$) & 0.915 & 0.946 & 0.910 & 0.625 & 0.308 &
0.203 \\
[2pt]
& Adaptive Neyman &0.948 & 0.960 & 0.967 & 0.952 & 0.843 & 0.613 \\
%
%
%
\hline
\end{tabular*}
\end{table}

Table \ref{tab:table2} shows that the empirical frequencies of finding
the true dimensions for different situations. As can be expected,
estimating the true dimension becomes more challenging as the model
becomes more complicated. For example, the probabilities of finding the
true dimension for model 3 are much smaller than those for model 2. Our
simulation results also show that, for a range of small values of $m$,
the two $\chi^2$ procedures perform very well especially if $n=500$,
where for brevity those results are represented by $m=5$ and $7$ in
Table \ref{tab:table2}. However, when $m=30$, the probabilities of
finding the true dimension become smaller for those procedures, which
is especially true for models 2 and 3. This is another illustration
that using a large number of principal components will lead to a loss
of power for the underlying $\chi^2$ tests. Again, the adaptive Neyman
procedure performs comparably to the two $\chi^2$ procedures with a
well-chosen $m$.


%

\subsection{Simulation 2: Sensitivity to the truncation point
$m$}\label{sec:simulation2} In the examples in Section
\ref{sec:simulation1}, the three EDR directions are linearly independent
when projected onto the three leading principal
components. Hence, the two $\chi^2$ procedures are expected to work
so long as $m\ge3$. For situations where one or more EDR directions only
depend on high-order principal components,
the choice of $m$ in the $\chi^2$ procedure is crucial and the
adaptive Neyman procedure has a
clear advantage.

To illustrate this, we consider a new model
\[
\mbox{Model 4: } Y={\langle}\beta_1, X{\rangle}\times(2{\langle
}\beta_4,X{\rangle}
+1)+\varepsilon,
\]
where $X$ is a Gaussian process whose distribution is as described in Section
\ref{sec:simulation1}, $\beta_1$ is as in (\ref{e:dir}), but
$\beta_4$ is given by
\begin{eqnarray*}
\beta_4(t)&=&0.45\sqrt{2}\cos(2\pi t)+0.6\sqrt{2}\cos(4\pi t)
-3\sqrt{2}\sin(6\pi t) \\
&&{} +1.2\sqrt{2}\sin(8\pi t)+\sum_{k>4} {(-1)^k \sqrt{2} \over
(2k)^3} \sin(2k\pi t).
\end{eqnarray*}
In this model, the dimension of the EDR space is 2, but the
projections of $\beta_1$ and $\beta_4$ onto the first five principal
components are linearly dependent; indeed, $K_{(m)}=1, m \le5$, and
$K_{(m)}=2, m\ge6$.
As shown in Table \ref{tab:model4},
the two $\chi^2$ procedures with $m=5$ both failed to find the true dimension,
even when $n=500$. On the other hand, when $n=500$ and $m=7$, the
$\chi^2$ procedures worked very well. With $m=30$, both $\chi^2$
tests again have considerably
lower powers, which leads to smaller probabilities of correct identification.
As in the previous models, the adaptive Neyman procedure has comparable
performance
to the best $\chi^2$ procedures.

%
\begin{table}
\tablewidth=240pt
\caption{Empirical frequencies of finding the
correct model in model 4} \label{tab:model4}
\begin{tabular*}{\tablewidth}{@{\extracolsep{\fill}}lcc@{}}
\hline
& $\bolds{n=200}$ & $\bolds{n=500}$\\
\hline
%
$\chi^2$ test ($m=5$) &0.040 & 0.047 \\
Adj. $\chi^2$ ($m=5$) &0.068 & 0.068 \\
[2pt]
$\chi^2$ test ($m=7$) &0.358 & 0.913 \\
Adj. $\chi^2$ ($m=7$) &0.410 & 0.899 \\
[2pt]
$\chi^2$ test ($m=30$) &0.085 & 0.566 \\
Adj. $\chi^2$ ($m=30$) &0.170 & 0.616 \\
[2pt]
Adaptive Neyman &0.229 & 0.885 \\
\hline
\end{tabular*}
\end{table}

%
\begin{figure}[b]

\includegraphics{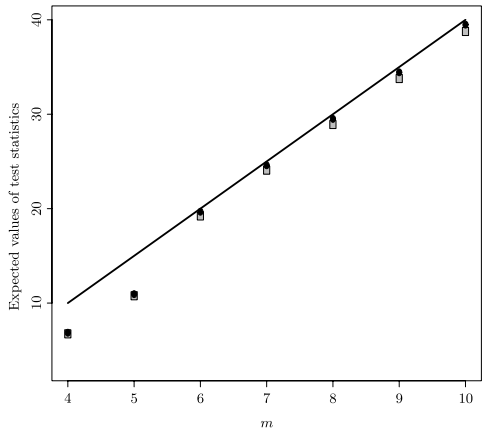}

\caption{The expected\vspace*{1pt} values of the test statistics $\CT_{K_0,(m)}$
and $\CT_{K_0,(m)}^\ast$ plotted as a function of $m$; the solid line
describes the theoretical
expected values, the rectangles are the means of $\CT_{K_0,(m)}$ and the
circles are the means of $\CT_{K_0,(m)}^\ast$.}
\label{fig:testmean}
\end{figure}

Finally, we use model 4 to illustrate the null distribution of $\CT
_{K_0,(m)}$ and $\CT_{K_0,(m)}^\ast$ when $K_{(m)} < K_0$. Consider
$K_0=2$; for each $m=4,5,\ldots,$ compute the expected values of
$\CT_{K_0,(m)}$ and $\CT_{K_0,(m)}^\ast$ by simulations and compare the
expectations with theoretical expectations $(m-K_0)(H-K_0-1)$. The
results are described in Figure~\ref{fig:testmean}, in which the grey
rectangles mark the means of $\CT_{K_0,(m)}$, the black circles mark
the means of $\CT_{K_0,(m)}^\ast$, and straight line represents
$(m-K_0)(H-K_0-1)$. The case $m\ge6$ correspond to (i) of\vspace*{1pt}
Theorem \ref{thm:fsir_asymp}, for which $\chi_{(m-K_0)(H-K_0-1)}^2$ is
the asymptotic distribution of the test statistics; the cases $m=4$ and
$5$ correspond to (ii) of Theorem \ref{thm:fsir_asymp}, for which $\chi
_{(m-K_0)(H-K_0-1)}^2$ is only a stochastic bound of the test
statistics. Both of these points are clearly reflected in Figure
\ref{fig:testmean}.


\subsection{Simulation 3: Multivariate data with measurement
errors}\label{s:multivariate}

In this subsection, we present a simulation study for high-dimensional
multivariate data; in particular,
we use the study to support claims made in (a) of Section \ref{s:discussions}.
Assume that $X$ is multivariate, and, for clarity, denote $X$ by $\BX$ here.
The simulated model is described as follows. We first generate a
$p$-dimensional variable
$\BX=(\BX_1^{\mathrm{T}},\BX_2^{\mathrm{T}})^{\mathrm{T}}$, where $p$ is ``large,'' with $\BX_1$
denoting a 10-dimensional random
vector while $\BX_2={\mathbf0}_{p-10}$, so that $\BX_1$ contains the
real signal in $\BX$.
Suppose that $\BX_1$ has a low-dimensional representation
\[
\BX_1=\sum_{k=1}^5 \xi_k \bolds{\psi}_k,
\]
where the $\bolds{\psi}_k$'s are orthonormal vectors, $\xi
_k\sim\Normal(0,\omega_k)$
are independent, and $(\omega_1,\ldots, \omega_5)=(3,2.8,2.6,2.4,
2.2)$; the
$\bolds{\psi}_k$'s are randomly generated, but are fixed
throughout the simulation study.
Furthermore, instead of observing the true~$\BX$, assume that we
observe an error-prone surrogate,
\[
\BW=\BX+\BU,
\]
where $\BU\sim\Normal({\mathbf0}, \mathbf{I}_p)$ is
measurement error.
Thus, the eigenvalues of the covariance of $\BW$ are bounded below by 1.
Note that this is a simpler measurement-error model than the one
considered in \citet{r7}, but
realistically portrays certain crucial aspects of high-dimensional data
encountered in practice; for example, in a typical fMRI study the total
number of brain voxels
captured by the image is huge but often only a relatively small portion
of the voxels are active for the task being studied, while background
noise is ubiquitous.

Let $X_1,\ldots,X_p$ be the components of $\BX$. Consider the model
\[
\mbox{Model 5: } Y=(X_1+X_2)/(X_2+X_3+X_4+X_5+1.5)^2+\varepsilon,
\]
where $\varepsilon\sim\Normal(0,0.5^2)$. Thus, $Y$ only depends on
$\BX_1$.
Below we compare the $\chi^2$ procedure in \citet{r22} and the adaptive
Neyman procedure using $\BW$ as the observed covariate.

%
\begin{table}[b]
\caption{Empirical frequencies of finding the
correct model in model 5} \label{tab:multivariate}
\begin{tabular*}{\tablewidth}{@{\extracolsep{\fill}}lccccc@{}}
\hline
& & $\bolds{p=15}$ & $\bolds{p=20}$ & $\bolds{p=40}$ & $\bolds{p=100}$
\\
\hline
$n=200$ & Li's $\chi^2$ test & 0.328 & 0.258 & 0.123 & 0.007 \\
& Adaptive Neyman & 0.588 & 0.596 & 0.562 & 0.528 \\
[2pt]
$n=500$ & Li's $\chi^2$ test & 0.898 & 0.833 & 0.612 & 0.276 \\
& Adaptive Neyman & 0.955 & 0.956 & 0.953 & 0.960 \\
\hline
\end{tabular*}
\end{table}

%
\begin{table}
\caption{Empirical sizes of Li's $\chi^2$ test and the adaptive
Neyman test for $H_0\dvtx K\le2$ under model 5}
\label{tab:multivariate_size}
\begin{tabular*}{\tablewidth}{@{\extracolsep{\fill}}lccccc@{}}
\hline
& & $\bolds{p=15}$ & $\bolds{p=20}$ & $\bolds{p=40}$ & $\bolds{p=100}$
\\
\hline
$n=200$ & Li's $\chi^2$ test & 0.015 & 0.012 & 0.002 & 0.001\\
& Adaptive Neyman & 0.016 & 0.013 & 0.009 & 0.010\\
[2pt]
$n=500$ & Li's $\chi^2$ test & 0.044 & 0.042 & 0.025 & 0.013 \\
& Adaptive Neyman & 0.037 & 0.035 & 0.035 & 0.030 \\
\hline
\end{tabular*}
\end{table}

We conducted simulations for $n=200, 500$ and $p=15, 20, 40, 100$. For
each setting, we repeat the simulation for 1000
times and used Li's procedure and the adaptive Neyman procedure in
deciding the number of EDR directions.
For both procedures, the nominal size $\alpha=0.05$ was used. In Table
\ref{tab:multivariate}, we summarize the
empirical frequencies of finding the correct dimension. As can be seen,
while the performance of the adaptive
Neyman procedure is quite stable for different $p$'s, the performance
of Li's procedure deteriorates
as $p$ increases. In Table \ref{tab:multivariate_size}, we also
present the true sizes, obtained by
simulations, of the two tests
for $H_0\dvtx K\le2$. In all cases both tests have sizes under 0.05.
The sizes of both tests are closer to the nominal size when $n=500$
than when $n=200$.
With a fixed $n$, the sizes of Li's test decrease quickly as $p$ increases,
reflecting the conservative nature of the test for large $p$,
while those for the adaptive Neyman test remain relatively stable.

%
\section{Data analysis}\label{sec:data}

In this section, we consider the Tecator data [Thodberg (\citeyear
{r31})], which
can be
downloaded at \href{http://lib.stat.cmu.edu/datasets/tecator}{http://lib.stat.cmu.edu/datasets/tecator}.
The data were previously analyzed in a number of papers including
\citet{r14},
\citet{r1} and Hsing and Ren (\citeyear{r19}).
The data contains measurements obtained by analyzing
215 meat samples,
where for each sample a 100-channel, near-infrared spectrum was
obtained by a
spectrometer, and the water, fat and protein contents were also directly
measured. The spectral data can be naturally
considered as functional data, and we are interested in building a
regression model to predict the fat content from the spectrum.
Following the convention in the literature, we applied a logistic
transformation to the percentage of fat
content, $U$, by letting $Y=\log_{10} \{U/(1-U)\}$.

In applying functional SIR, both \citet{r14} and \citet{r1}
used graphical tools to select
the number EDR directions, where the numbers of directions selected
were 10 and 8, respectively. On the other hand, using only two EDR directions,
\citet{r1} applied MAVE to achieve a
prediction error comparable to what can be achieved by SIR using 8
directions. These
conclusions were somewhat inconsistent.

\begin{figure}

\includegraphics{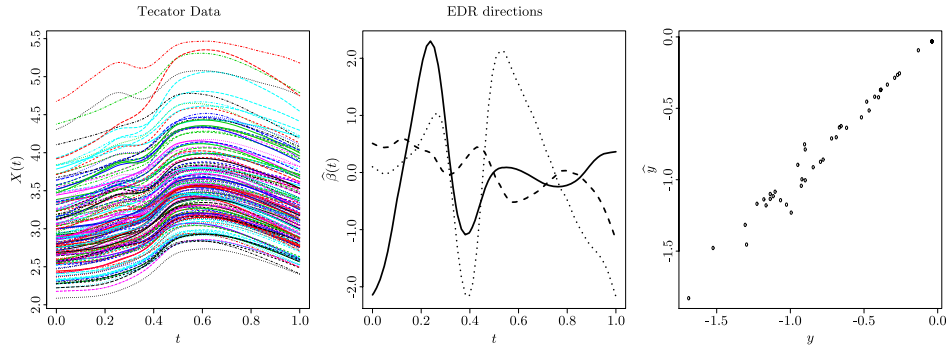}

\caption{Tecator spectrum data: the first plot shows the
Tecator spectrum data in the training set, the second plot show the estimated
EDR directions and the last plot is the
predicted vs. true fat contents for the test data set.}
\label{fig:tecator}
\end{figure}

Based on the instructions given by the Tecator website, we used the
first 172 samples for training and the last $43$ for testing.
Following \citet{r1}, we focused on the most informative
part of the spectra, with wavelengths ranging from 902 to 1028 nm.
The curves are rescaled onto the interval $[0,1]$. The first plot in
Figure \ref{fig:tecator}
shows those spectra in the training set.

%

We first fitted B-splines to the discrete data, and then applied
our\break
sequential-testing procedures.
With $\alpha=0.05$, the adaptive Neyman procedure concluded that $K=3$.
To see how well a three-dimensional model works, the model
$Y=f({\langle}
\beta_1, X{\rangle},{\langle}\beta_2,
X{\rangle}, {\langle}\beta_3, X{\rangle})+\varepsilon$ was
entertained. The EDR
directions were estimated by the
regularized approach, RSIR, introduced by \citet{r37};
the estimated EDR directions are presented in the center
plot in Figure \ref{fig:tecator}. Finally, the link function $f$ was
estimated by smoothing spline ANOVA
[\citet{r15}] with interaction terms; the estimated model was then applied
to test data
to predict fat content. The root mean prediction error was 0.062
which is comparable to what was obtained by MAVE in \citet{r1}.
The plot of the predicted versus the true fat contents for test
data is also given in
Figure \ref{fig:tecator}.

Our result is in agreement with what was obtained using MAVE in
\citet{r1}
in that a low-dimensional model is appropriate for this data set.
\begin{appendix}\label{app}
\section*{Appendix: Proofs}

In the following, the notation ``$\star$'' refers to symbolic matrix
multiplication;
for instance, if $f_1,\ldots,f_k$ are mathematical objects (functions,
matrices, etc.)
and $\bdc=(c_1,\ldots,c_k)\tran$ is a vector for which
the operation $\sum_{i=1}^k c_if_i$ is defined, we will denote the sum
by $(f_1,\ldots,f_k)\star\bdc$;
also if $C$ is a matrix containing columns $\bdc_1, \ldots, \bdc
_\ell$, the notation
$(f_1,\ldots,f_k)\star C$ refers the array
$[(f_1,\ldots,f_k)\star\bdc_1,\ldots,(f_1,\ldots,f_k)\star\bdc
_\ell]$.

The notation is defined in Section \ref{sec:fsir} and will be used
extensively below
without further mention.

\subsection{\texorpdfstring{Proof of Theorem
\protect\ref{thm:fsir_asymp}}{Proof of Theorem 3.1}}
Recall from (\ref{e:VB}) that $V_{(m)} = B_{(m)} B_{(m)}\tran$ and
$\widehat V_{(m)} = \widehat B_{(m)}\widehat B_{(m)}\tran$.
Thus, to study the eigenvalues of $V_{(m)}$ and $\widehat V_{(m)}$,
we can equivalently study the singular values of $B_{(m)}$ and $
\widehat
B_{(m)}$, respectively.
Recall that $K_{(m)}=\operatorname{Rank}(B_{(m)})$.
Under the hypothesis $K\le K_0$, we also have $K_{(m)}\le K_0$. Suppose
$B_{(m)}$ has the following singular-value decomposition:
\[
B_{(m)}=\CP\pmatrix{
D &0\cr0 &0
} \CQ\tran,
\]
where $D:=\diag( \lambda_1^{1/2}(V_{(m)}), \ldots,
\lambda_{K_{(m)}}^{1/2}(V_{(m)}))$ contains
the nonzero singular values of $B_{(m)}$, and $\CP$ and $\CQ$ are
orthonormal matrices
of dimensions $m\times m$ and $H\times H$, respectively, which contain
the singular
vectors of $B_{(m)}$. Note that, for brevity of notation, we leave out
$m$ in $\CP,\CQ$ and $n$
in $\widehat B_{(m)}, \widehat V_{(m)}$ in this proof.

Partition $\CP$ and $\CQ$ as $\CP= [ \CP_1 | \CP_2 ]$, $\CQ=
[ \CQ_1 | \CQ_2 ]$
where $\CP_1$ and $\CQ_1$ both have $K_{(m)}$ columns, and
$\CP_2$ and $\CQ_2$ have $m-K_{(m)}$ and $H-K_{(m)}$ columns, respectively.
Thus, the columns of $\CP_2$ and $\CQ_2$ are singular vectors
corresponding to the singular value
$0$, and so $B_{(m)}\tran\CP_2= \mathbf{0}$ and
$B_{(m)}\CQ_2=\mathbf{0}$. We further partition $\CQ_2$ in the
following way.
Recall that $\mu_1, \ldots,\mu_H$ are the within-slice means
defined in (\ref{eq:within_slice_mean_curve}). By Theorem \ref{thm:fsir},
$\operatorname{span}(\mu_1-\mu,\ldots,\mu_H-\mu)$ is a subspace of
$\operatorname{span}(\Gamma_{X} \beta_1,\ldots, \Gamma_{X} \beta_K)$ and
therefore has dimension
less than or equal to $K\le K_0$. It follows from the ``rank-nullity
theorem'' that there exists a matrix
$\CQ_{2\circ}$ of dimension $H\times(H-K_0)$ with orthonormal
columns such that
%
%
\begin{equation}\label{eq:q2}
(\mu_1-\mu,\ldots,\mu_H-\mu) \star(F \CQ_{2\circ})= \mathbf{0}.
\end{equation}
Furthermore, observe that $\bdg$ spans the null space of $F$ and so
$(\mu_1-\mu,\ldots,\mu_H-\mu)\star(F\bdg) = 0$. Without loss of
generality, let
$\bdg$ be the last column of $\CQ_{2\circ}$.
Define an operator $T\dvtx\CH\to\bbR^m$ by
%
%
\begin{equation}\label{e:T}
T x= (\omega_1^{-1/2} {\langle}\psi_1, x{\rangle}, \ldots, \omega_m^{-1/2}
{\langle}\psi_m, x{\rangle})\tran,\qquad x\in\CH.
\end{equation}
Applying $T$ to both sides of (\ref{eq:q2}), we have
%
%
\begin{equation}\label{e:BQ2}
B_{(m)}\CQ_{2\circ}=MF\CQ_{2\circ}=T\{(\mu_1-\mu,\ldots,\mu
_H-\mu)\} F \CQ_{2\circ}
= \mathbf{0},
\end{equation}
where, for convenience, the notation $T\{(\mu_1-\mu,\ldots,\mu
_H-\mu)\}$ means
$(T(\mu_1-\mu),\ldots,T(\mu_H-\mu))$.
This means $\CQ_{2\circ}$ is contained in the column space of $\CQ_{2}$.
Without loss of generality, we assume that $\CQ_2$ has the decomposition
%
%
\begin{equation}\label{e:partitionQ}
\CQ_2 = [ \CQ_{2*} | \CQ_{2\circ} ],
\end{equation}
where $\CQ_{2*}$ is of dimension $H\times(K_0-K_{(m)})$.
When $m$ is large enough so that $K_{(m)}=K_0$, then $\CQ_2=\CQ
_{2\circ}$.

%
%
%
%
%

Let
%
%
\begin{equation}\label{e:U_n}
U_{n} := \CP_2\tran\widehat B_{(m)} \CQ_2=\CP_2\tran\bigl( \widehat
B_{(m)}-B_{(m)}\bigr) \CQ_2.
\end{equation}
Also define
%
%
\begin{equation}\label{e:notation2}
\widetilde\bdvartheta_{h,(m)} = {1\over n_h}\sum_i \bdeta_{i,(m)}
I(Y_i\in
S_h) \quad\mbox{and}\quad
\widetilde M = \bigl[\widetilde\bdvartheta_{1,(m)}, \ldots,\widetilde
\bdvartheta
_{H,(m)}\bigr]_{m\times H}.\hspace*{-30pt}
\end{equation}
\begin{lm}\label{lm:asymp_m}
Assume that $X(t)$ has an elliptically contoured distribution
satisfying (\ref{e:sph}).
Let $\widetilde M$ be defined by (\ref{e:notation2}). We have
$\sqrt{n} \CP_2\tran\widetilde M G\cid\CZ\Lambda$, where $\CZ$
is a
$(m-K_{(m)})\times H$ matrix of independent $\Normal(0,1)$ random
variables, $\Lambda=\diag(\tau_1^{1/2},\ldots,\tau_H^{1/2})$ with
$\tau_h=\E(\Theta^2|Y\in S_h)$.
\end{lm}
\begin{pf} Let
\[
\bdu_{n,h} = {1\over n}\sum_{i=1}^n \bdeta_{i,(m)} I(Y_i\in S_h),\qquad
p_{n,h} = {1\over n}\sum_{i=1}^n I(Y_i\in S_h),
\]
then $\widetilde\bdvartheta_{h,(m)} = \bdu_{n,h} /p_{n,h}$, denote
$\bdu
_h=\E(\bdu_{n,h})=\bdvartheta_{h,(m)} p_h$.
Then
\[
\widetilde M = \biggl({\bdu_{n,1}\over p_{n,1}},\ldots,{\bdu
_{n,H}\over
p_{n,H}} \biggr)_{m\times H}
\]
and
\[
M = \biggl({\bdu_{1}\over p_{1}},\ldots,
{\bdu_{H}\over p_{H}} \biggr)_{m\times H},
\]
and so
%
%
\begin{eqnarray}\label{e:Z1}
&&
n^{1/2}(\widetilde M-M) \nonumber\\
&&\qquad= n^{1/2} \biggl({\bdu_{n,1}\over
p_{n,1}}-{\bdu
_{1}\over p_{1}},\ldots,
{\bdu_{n,H}\over p_{n,H}}-{\bdu_{H}\over p_{H}} \biggr) \nonumber\\
&&\qquad= n^{1/2} \biggl({\bdu_{n,1}-\bdu_1\over p_{n,1}},\ldots,
{\bdu_{n,H}-\bdu_H\over p_{n,H}} \biggr) \nonumber\\
&&\qquad\quad{} -n^{1/2} \biggl({\bdu_{1}\over p_{n,1}p_1}(p_{n,1}-p_1),\ldots,
{\bdu_{H}\over p_{n,H}p_H}(p_{n,H}-p_H) \biggr) \\
&&\qquad=n^{1/2} \biggl({\bdu_{n,1}-\bdu_1\over p_{1}},\ldots,
{\bdu_{n,H}-\bdu_H\over p_{H}} \biggr) \nonumber\\
&&\qquad\quad{}-n^{1/2} \biggl({\bdu_1\over p_1^2}(p_{n,1}-p_1),\ldots,
{\bdu_H\over p_H^2}(p_{n,H}-p_H) \biggr)\nonumber\\
&&\qquad\quad{} + o_p(1). \nonumber
\end{eqnarray}
By the central limit theorem and covariance computations, it is easy to
see that
the columns of $\sqrt{n} \CP_2\tran(\widetilde M-M) G$ are
asymptotically independent, where the $h$th column converges in
distribution to a random vector having the multivariate normal distribution
%
%
\begin{equation}\label{e:normal}
\Normal\bigl(\mathbf{0}, \var\bigl(\CP_2\tran\bdeta_{(m)}|Y\in S_h\bigr) \bigr).
\end{equation}
For convenience, let $\CV$ denote the vector $({\langle}\beta
_1,X{\rangle},
\ldots, {\langle}\beta_K,X{\rangle})$.
By iterative conditioning,
\begin{eqnarray*}
\var\bigl(\CP_2\tran\bdeta_{(m)}|Y\in S_h\bigr)
&=& \E\bigl\{\var\bigl(\CP_2\tran\bdeta_{(m)}|\CV, \Theta, Y,\varepsilon
\bigr)|Y\in S_h\bigr\} \\
&&{}+ \var\bigl\{\CP_2\tran\E\bigl(\bdeta_{(m)}|\CV,\Theta, Y,\varepsilon
\bigr)|Y\in S_h\bigr\} \\
&=& \E\bigl\{\var\bigl(\CP_2\tran\bdeta_{(m)}|\CV, \Theta\bigr)|Y\in S_h\bigr\} \\
&&{}+ \var\bigl\{\CP_2\tran\E\bigl(\bdeta_{(m)}|\CV,\Theta\bigr)|Y\in S_h\bigr\},
\end{eqnarray*}
where we used the facts that $Y$ is redundant given $\varepsilon$ and
the ${\langle}\beta_k,X{\rangle}$'s,
and $X$ is independent of $\varepsilon$.
With the notation $\check\CV=({\langle}\beta_1,\check X{\rangle},
\ldots, {\langle}
\beta_K,\check X{\rangle})$, we have
%
%
\begin{eqnarray}\label{e:vardecomp}
\var\bigl(\CP_2\tran\bdeta_{(m)}|Y\in S_h\bigr)
&=& \E\bigl\{\Theta^2\var\bigl(\CP_2\tran\check\bdeta_{(m)}|\check\CV
\bigr)|Y\in S_h\bigr\} \nonumber\\[-8pt]\\[-8pt]
&&{} + \var\bigl\{\Theta^2\E\bigl(\CP_2\tran\check\bdeta_{(m)}|\check\CV
\bigr)|Y\in S_h\bigr\}. \nonumber
\end{eqnarray}
In the following, we focus on the special case $\check X$ is Gaussian.
The general case
is similar but requires a more careful analysis of the conditional
distribution of
jointly elliptically contoured random variables. Let $\bdb$ be any
column of $\CP_2$. Then
\[
\E\bigl(\bdb\tran\check\bdeta_{(m)}\bigr) = 0 \quad\mbox{and}\quad
\E\bigl(\bdb\tran\check\bdeta_{(m)}{\langle}\beta_k,\check X{\rangle
}\bigr) = \bdb
\tran\bdb_{k,(m)} = 0,\qquad 1\le k\le K.
\]
Thus, $\CP_2\tran\check\bdeta_{(m)}$ is a vector of standard normal
random variables that
are independent of the ${\langle}\beta_k,\check X{\rangle}$'s. It
follows from
(\ref{e:vardecomp}) that
%
%
\begin{equation}\label{e:vardecomp1}
\var\bigl(\CP_2\tran\bdeta_{(m)}|Y\in S_h\bigr)
= \E(\Theta^2|Y\in S_h) I = \tau_h I,
\end{equation}
where $I$ is the identity matrix. The proof is complete.
\end{pf}
\begin{lm}\label{lm:approx}
Let $X(t), \CZ$ and $\Lambda$ be as in Lemma \ref{lm:asymp_m}. Then
$\sqrt{n} U_{n}\cid Z$ where all of the entries of $Z$ are normally
distributed with mean zero
and have the following properties:
\begin{longlist}
\item
If $K_{(m)}=K_0$, then $Z\eid\CZ\Lambda\CJ_{\bdg}\CQ_{2}$.
\item
If $K_{(m)} < K_0$, then $Z$ can be partitioned as $Z=[ Z_* |
Z_{\circ} ]$ in accordance
with the partition of $\CQ_2$ in (\ref{e:partitionQ}), where
$Z_{\circ}\eid\CZ\Lambda\CJ_{\bdg}\CQ_{2\circ}$;
furthermore, if $X$ is Gaussian then $Z_*$ and $Z_{\circ}$ are
independent, where the last
column of $Z_{\circ}$ is identically 0 while the rest of the entries
of $Z_{\circ}$ are i.i.d. standard normal.
\end{longlist}
\end{lm}
\begin{pf} First, write 
\begin{eqnarray*}
n^{1/2} U_{n}&=& n^{1/2} \CP_2\tran\widehat M \widehat F \CQ_2\\
&=& n^{1/2} \CP_2\tran\{\widetilde M F+(\widehat M-\widetilde M) F +
(\widehat M-\widetilde M)(\widehat
F-F)+ \widetilde M(\widehat F-F)\} \CQ_2.
\end{eqnarray*}
Denote $\bar X_h= \sum_i X_i I(Y_i\in S_h)/\sum_i I(Y_i\in S_h)$.
Then
\[
\widehat M = \{ {\langle}\widehat\omega_j^{-1/2}\widehat\psi
_j,\bar X_h-\bar X{\rangle}\}
_{j,h=1}^{m,H},\qquad
\widetilde M = \{ {\langle}\omega_j^{-1/2}\psi_j,\bar X_h-\mu
{\rangle}\}_{j,h=1}^{m,H}.
\]
Then
%
%
\begin{eqnarray}\label{e:Mdiff}
\widehat M-\widetilde M
&=& \{ {\langle}\widehat\omega_j^{-1/2}\widehat\psi_j-\omega
_j^{-1/2}\psi_j,\bar
X_h-\bar X{\rangle}\}_{j,h=1}^{m,H}\nonumber\\[-8pt]\\[-8pt]
&&{} - \{ {\langle}\omega_j^{-1/2}\psi_j,\bar X-\mu{\rangle}\}_{j,h=1}^{m,H}.
\nonumber
\end{eqnarray}
It follows that
%
%
\begin{eqnarray}\label{e:hall}
\widehat\psi_j(t)-\psi_j(t)&=& \sum_{\ell\neq j} {\psi_\ell(t)
\over
\omega_j-\omega_\ell}
{\langle}(\widehat\Gamma_X-\Gamma_X) \psi_\ell,\psi_j {\rangle}
+O_p(n^{-1}),\nonumber\\[-8pt]\\[-8pt]
\widehat\omega_j-\omega_j&=& {\langle}(\widehat\Gamma_X-\Gamma
_X) \psi_\ell
,\psi_j {\rangle}
+O_p(n^{-1}). \nonumber
\end{eqnarray}
These were established by (2.8) and (2.9) in \citet{r16} for $\CH
=L^2[a,b]$.
Actually, they hold for any Hilbert space $\CH$; see \citet{r11a},
Theorem 3.8.11.
Since $\widehat\Gamma_X-\Gamma_X = O_p(n^{-1/2})$, these imply
$\widehat\omega_j=\omega_j+O_p(n^{-1/2})$ and\vspace*{2pt} $\widehat\psi_j=\psi
_j+O_p(n^{-1/2})$.
Also $\bar X-\mu= O_p(n^{-1/2})$. Thus, $\widehat M-\widetilde M =
O_p(n^{-1/2})$.
Since we also have $\widehat F-F=O_p(n^{-1/2})$, we conclude that
\[
n^{1/2} \CP_2\tran(\widehat M-\widetilde M) (\widehat F-F) = O_p(n^{-1/2}).
\]
Similarly, since $\CP_2\tran M = 0$,
\[
n^{1/2} \CP_2\tran\widetilde M (\widehat F-F) = n^{1/2} \CP_2\tran
(\widetilde M -
M)(\widehat F-F) = O_p(n^{-1/2}).
\]
Thus,
%
%
\begin{equation}\label{e:Un}
n^{1/2} U_{n}= n^{1/2} \CP_2\tran\widetilde M F \CQ_2 +n^{1/2}\CP_2
(\widehat
M-\widetilde M) F \CQ_2
+o_p(1).
\end{equation}
To get the desired result, we break the proof into several parts.
First, we establish that
%
%
\begin{equation}\label{e:asympA}\quad
\bigl(n^{1/2} \CP_2\tran(\widetilde M-M) F \CQ_2, n^{1/2}\CP_2
(\widehat
M-\widetilde M) F \CQ_2 \bigr)
\cid(Z_1,Z_2),
\end{equation}
where $(Z_1,Z_2)$ are jointly normal with mean $0$.
By (\ref{e:Mdiff}) and the fact that $(1,\ldots,1)F=\mathbf{0}$,
we have
\[
n^{1/2}(\widehat M-\widetilde M)F \CQ_2
= n^{1/2} \{ {\langle}\widehat\omega_j^{-1/2}\widehat\psi_j-\omega
_j^{-1/2}\psi
_j,\bar X_h-\bar X{\rangle}\}_{j,h=1}^{m,H}
F\CQ_2.
\]
Since $\widehat\omega_i^{-1/2}\widehat\psi_i-\omega_i^{-1/2} \psi
_i=O_p(n^{-1/2})$
and $\bar X_h-\bar X \cip\mu_h-\mu$,
%
%
\begin{eqnarray}\label{e:asympC}\quad
n^{1/2}(\widehat M-\widetilde M)F \CQ_2
&=& n^{1/2} \{ {\langle}\widehat\omega_j^{-1/2}\widehat\psi_j-\omega
_j^{-1/2}\psi
_j,\gamma_h{\rangle}\}_{j,h=1}^{m,H-K_{(m)}}\nonumber\\[-8pt]\\[-8pt]
&&{} + o_p(1),\nonumber
\end{eqnarray}
where
\[
(\gamma_1,\ldots,\gamma_{H-K_{(m)}}) = (\mu_1-\mu,\ldots,\mu
_H-\mu)\star(F\CQ_2).
\]
Note that the last $H-K_0$ of the $\gamma_k$'s are equal to 0 by (\ref
{eq:q2}).
In particular, if $K_{(m)}=K_0$ then all of the $\gamma_h$'s are equal
to 0 and
(\ref{e:asympA}) is established with $Z_2=0$ and
$Z_1\eid\CZ\Lambda\CT_{\bdg}\CQ_{2}$ by Lemma \ref{lm:asymp_m}.
The assertion (i) follows readily from (\ref{e:Un}). Below, we focus
on the case $K_{(m)}<K_0$.
Recall that
\[
\{ {\langle}\omega_j^{-1/2}\psi_j,\gamma_h{\rangle}\}
_{j,h=1}^{m,H-K_{(m)}}=
MF\CQ_2 = 0,
\]
which implies that
%
%
\begin{equation}\label{e:gamma}
\{ {\langle}\psi_j,\gamma_h{\rangle}\}_{j,h=1}^{m,H-K_{(m)}} = 0.
\end{equation}
By (\ref{e:asympC}) and (\ref{e:gamma}),
%
%
\begin{equation}\label{e:asympB}
n^{1/2}(\widehat M-\widetilde M)F \CQ_{2*}
= n^{1/2} \{ {\langle}\omega_j^{-1/2}(\widehat\psi_j-\psi
_j),\gamma_h{\rangle}\}
_{j,h=1}^{m,K_0-K_{(m)}}
+ o_p(1).\hspace*{-34pt}
\end{equation}
By the central limit theorem, the random element
$n^{-1/2}(\widehat R-R,\bdu_{n,h}-\bdu_h,
p_{n,h}-p_h, h = 1,\ldots,H)$ has a jointly Gaussian limit.
In view of (\ref{e:Z1}), (\ref{e:hall}) and (\ref{e:asympB}), the
claim in
(\ref{e:asympA}) is established by performing a linear transformation.
Define the partitions
$Z_1 = [ Z_{1*} | Z_{1\circ} ]$ and $Z_2 = [ Z_{2*} | Z_{2\circ
} ]$ and so
$[ Z_* | Z_{\circ} ] = [ Z_{1*}+Z_{2*} | Z_{1\circ} ]$ since
$Z_{2\circ}= 0$.
By Lemma \ref{lm:asymp_m},
%
%
\begin{equation}\label{e:asympD}
[ Z_{1*} | Z_{1\circ} ] \eid[ \CZ\Lambda\CJ_{\bdg}\CQ_{2*}
| \CZ\Lambda\CT_{\bdg}\CQ_{2\circ} ]
\end{equation}
and so $Z_{\circ}=Z_{1\circ}\eid\CZ\Lambda\CJ_{\bdg}\CQ_{2\circ}$.
Assume for the rest of the proof that $X$ is Gaussian.
Recall that $\CJ_{\bdg}=I-\bdg\bdg\tran$. By the fact that $\tau
_h\equiv1$ and
the convention that $\bdg$ is the last column of $\CQ_2$, it follows
from (\ref{e:asympD}) that
\[
[ Z_{1*} | Z_{1\circ} ] \eid[ \CZ\CQ_{2*} | \CZ\widetilde
\CQ
_{2\circ} ],
\]
where $\widetilde\CQ_{2\circ}$ denotes the matrix whose last column
contains $0$'s but the
remaining entries are taken after $\CQ_{2\circ}$. By Lemma \ref
{lm:asymp_m}, $Z_{1*}$ and $Z_{1\circ}$ are
independent. So it remains to show that $Z_{2*}$ and $Z_{1\circ}$ are
independent. By (\ref{e:hall}) and (\ref{e:gamma}), with $\gamma
_{k\ell} := {\langle}\gamma_k,\psi_{\ell}{\rangle}$,
%
%
\begin{eqnarray}\label{e:asympD1}
&& n^{1/2}{\langle}\omega_j^{-1/2}(\widehat\psi_j-\psi_j),\gamma
_k{\rangle}
\nonumber\\
&&\qquad= n^{1/2}\omega_j^{-1/2} \sum_{\ell=m+1}^\infty{\gamma_{k\ell}
\over
\omega_j-\omega_\ell} \int(\widehat R-R) \psi_\ell\psi_j + o_p(1)
\nonumber\\
&&\qquad= n^{1/2} \omega_j^{-1/2} \sum_{\ell=m+1}^\infty{\gamma_{k\ell}
\over\omega_j-\omega_\ell}
\Biggl\{{1 \over n} \sum_{i=1}^n \xi_{i\ell}\xi_{ij} \Biggr\} +
o_p(1) \\
&&\qquad= n^{-1/2} \sum_{i=1}^n \Biggl(\sum_{\ell=m+1}^\infty{\gamma
_{k\ell}
\over\omega_j-\omega_\ell} \xi_{i\ell} \Biggr) \eta_{ij}+o_p(1)
\nonumber\\
&&\qquad=: z_{jk} + o_p(1), \nonumber
\end{eqnarray}
for $j=1,\ldots,m$, $k=1,\ldots,K_0-K_{(m)}$. Let $\bdz
_{k}=(z_{1k},\ldots,z_{mk})\tran$.
By (\ref{e:asympB}) and (\ref{e:asympD1}),
%
%
\begin{equation}\label{e:2nd}
n^{1/2} \CP_2\tran(\widehat M-\widetilde M)F \CQ_{2*} =
\CP_2\tran\bigl[\bdz_{1}, \ldots, \bdz_{K_0-K_{(m)}}\bigr] + o_p(1).
\end{equation}
Since $MF\CQ_2=0$ and $\CP_2^T\bdu_h=0$, it follows from (\ref
{e:Z1}) that
%
%
\begin{equation}\label{e:1st}
n^{1/2} \CP_2\tran(\widetilde M-M)F \CQ_{2\circ}
= n^{1/2}\CP_2\tran\biggl({\bdu_{n,1}\over p_{1}},\ldots,
{\bdu_{n,H}\over p_{H}} \biggr)F \CQ_{2\circ}
+ o_p(1).\hspace*{-28pt} 
\end{equation}
We now compute the covariances between the components of
(\ref{e:2nd}) and (\ref{e:1st}). Since the components are jointly normal,
our goal is to show that the covariances are all $0$.
Let $\bdq$ be a column of $\CQ_{2\circ}$ and $k=1,\ldots,K_0-K_{(m)}$.
Note that
\[
%
\CP_2\tran\bdz_{k}={1\over n^{1/2}} \sum_{i=1}^n \CP_2\tran\CD
_{ik} \bdeta_{i,(m)} ,
\]
where
\[
\CD_{ik}=\diag\Biggl(\sum_{\ell=m+1}^\infty{\gamma_{\ell k}
\over\omega_j-\omega_\ell} \xi_{i\ell}, j=1,\ldots,m \Biggr).
\]
Since $\E(\CP_2\tran\bdz_{k})=0$,
\begin{eqnarray*}
&& \cov\biggl(n^{1/2}\CP_2\tran\biggl({\bdu_{n,1}\over p_{1}},\ldots,
{\bdu_{n,H}\over p_{H}} \biggr)F\bdq,\CP_2\tran\bdz_{k} \biggr) \\
&&\qquad= \E\biggl(n^{1/2}\CP_2\tran\biggl({\bdu_{n,1}\over p_{1}},\ldots,
{\bdu_{n,H}\over p_{H}} \biggr)F\bdq\bdz_{k}\tran\CP_2 \biggr) \\
&&\qquad= \E\biggl(n^{1/2} \biggl({\CP_2\tran\bdu_{n,1}\bdz_{k}\tran\CP
_2\over p_{1}},\ldots,
{\CP_2\tran\bdu_{n,H}\bdz_{k}\tran\CP_2\over p_{H}} \biggr)\star
(F\bdq) \biggr).
\end{eqnarray*}
Let $\CV$ be as defined in the proof of Lemma \ref{lm:asymp_m}. By
the same conditioning argument
employed there,
%
%
\begin{equation}\label{e:iterexp}\quad
\E( n^{1/2}\CP_2\tran\bdu_{n,h} \bdz_{k}\tran\CP_2)
= \E\bigl[\E\bigl\{ \CP_2\tran\bdeta_{(m)} \bdeta_{(m)}\tran
\CD_{k} \CP_2 | \CV\bigr\}I(Y\in S_h)\bigr],
\end{equation}
where the sub-index $i$ is suppressed from the symbols on the
right-hand side since it
suffices to deal with a generic process $(X,Y)$ in computing expectations.
Note that $\CP_2\tran\bdeta_{(m)}$ and $\CV$ are independent,
$\bdeta_{(m)}$ and $\CD_{k}$ are independent, and $\bdeta_{(m)},
\CV, \CD_{k}$ are normally distributed with mean zero.
Then it is easy to conclude from (\ref{e:iterexp}) that
%
%
\begin{eqnarray}\label{e:iterexp1}
&&
\E( n^{1/2}\CP_2\tran\bdu_{n,h} \bdz_{k}\tran\CP_2)\nonumber\\
&&\qquad= \E\bigl(\CP_2\tran\bdeta_{(m)}\bdeta_{(m)}\tran\CP_2\bigr) \E[\E\{
\CP_2
\CD_{k} \CP_2 | \CV\}I(Y\in S_h)] \\
&&\qquad= \E[\E\{\CP_2\CD_{k} \CP_2 | \CV\}I(Y\in S_h)]. \nonumber
\end{eqnarray}
%
%
%
%
%
By the property of the normal distribution, each (diagonal) element of
$\E\{\CD_{k}|\CV\}$ can be written as
$\sum_{j=1}^K c_{j}{\langle}\beta_j, X-\mu{\rangle}$ for some
$c_{j},1\le j\le K$.
For convenience, denote $\E\{\CD_{k}|\CV\}$ as $T(X-\mu)$ where $T$
is a linear functional.
Thus,
%
%
\begin{equation}\label{e:iterexp2}
\E[\E\{ \CP_2\tran\CD_{k} \CP_2 | \CV\}I(Y\in S_h)]
= p_h \CP_2\tran T(\mu_h-\mu)\CP_2.
\end{equation}
As a result,
\begin{eqnarray*}
&& \E\biggl(n^{1/2} \biggl({\CP_2\tran\bdu_{n,1}\bdz_{k}\tran\CP
_2\over p_{1}},\ldots,
{\CP_2\tran\bdu_{n,H}\bdz_{k}\tran\CP_2\over p_{H}} \biggr)\star
(F\bdq) \biggr) \\
&&\qquad= \bigl(\CP_2\tran T(\mu_1-\mu)\CP_2,\ldots,\CP_2\tran T(\mu_H-\mu
)\CP_2\bigr)\star(F\bdq) \\
&&\qquad= \CP_2\tran\bigl(\bigl(T(\mu_1-\mu),\ldots, T(\mu_H-\mu)\bigr)\star
(F\bdq) \bigr) \CP_2
= 0,
\end{eqnarray*}
by (\ref{eq:q2}).
This shows that the covariances between the components of
(\ref{e:2nd}) and (\ref{e:1st}) are all equal to $0$, and
concludes the proof that $Z_{1\circ}$ and $Z_{2*}$ are independent.
\end{pf}
\begin{pf*}{Proof of Proposition \ref{pp:eign_value_bound}}
Assume for convenience that $Z_1$ has full column rank. If this is not
the case, a slight modification of
the proof below suffices.
Denote the $j$th column of $Z$ as $\bdz_j$, and construct
orthonormal vectors by applying the Gram--Schmidt
orthonormalization to the columns of $Z$:
\[
\bdv_1={\bdz_1\over\|\bdz_1\|},\qquad \bdv_j={ (I-\Pi_{j-1})\bdz
_j \over\|(I-\Pi_{j-1})\bdz_j\|},\qquad
j=2,\ldots, \min(p,q),
\]
where $\Pi_{j-1}=[\bdv_1,\ldots, \bdv_{j-1}] [\bdv_1,\ldots, \bdv
_{j-1}]\tran$ is the projection matrix to the space spanned by $\bdz
_1, \ldots, \bdz_{j-1}$. The following properties can be verified:
\begin{enumerate}[(a)]
\item[(a)] $\bdv_j\tran\bdz_k=0$ for all pairs $k<j$. This is the result
of the construction of the $\bdv_j$'s.
\item[(b)] 
$\bdv_j$ is independent of $\bdz_k$ for $k>\max(j,r)$. This follows
from the assumption on $Z_2$.
\item[(c)] 
$\bdv_j\tran\bdz_k \sim\Normal(0,1)$ for $k>\max(j,r)$. The proof
of this is easy: by (b)
and the fact that $\|\bdv_j\|=1$, $(\bdv_j^{\mathrm{T}}\bdz_k|\bdv_j) \sim
\Normal(0,1)$; since this conditional distribution
does not depend on $\bdv_j$, it is also the marginal distribution.

\item[(d)] 
$\bdv_j\tran\bdz_k$ and $\bdv_{j'}\tran\bdz_{k'}$ are
independent if
$k>\max(j,r)$ and $k'>\max(j',r)$.
The proof is as follows. First for the case $j, j', k < k'$, we have
\begin{eqnarray*}
\bbP(\bdv_j\tran\bdz_k\le x, \bdv_{j'}\tran\bdz_{k'}\le y)
&=& \E[\bbP(\bdv_j\tran\bdz_k\le x, \bdv_{j'}\tran\bdz
_{k'}\le y|\bdv_{j},\bdv_{j'},\bdz_k)] \\
&=& \E[I(\bdv_j\tran\bdz_k\le x)\bbP(\bdv_{j'}\tran\bdz
_{k'}\le y|\bdv_{j'})] \\
&=& \Phi(x)\Phi(y),
\end{eqnarray*}
where the last step follows from (c).
Next for $j, j'< k = k'$ and $j\not=j'$, we have
\[
\bbP(\bdv_j\tran\bdz_k\le x, \bdv_{j'}\tran\bdz_{k}\le y)
= \E[\bbP(\bdv_j\tran\bdz_k\le x, \bdv_{j'}\tran\bdz_{k}\le
y|\bdv_{j},\bdv_{j'})]
= \Phi(x)\Phi(y)
\]
since $\bdv_j\tran\bdv_{j'}=0$.
\item[(e)] $\bdv_j\tran\bdz_j=\|(I-\Pi_{j-1})\bdz_j\|$ for $j\ge r+1$
is the square root of a $\chi^2_{p-j+1}$ variable,
and it is independent of any $\bdv_{j'}\tran\bdz_{k'}$ with $k'>
j'\ge r$. The claims can be easily
verified using conditioning arguments similar to those in (c) and
(d).
\item[(f)] $\bdv_j\tran\bdz_j$ is independent of $\bdv_{j'}\tran\bdz
_{j'}$, for $j,j'\ge r+1$ and $j\neq j'$.
This can be verified by checking the independence between $(I-\Pi
_{j-1})\bdz_j$ and $(I-\Pi_{j'-1})\bdz_{j'}$.
\end{enumerate}
Based on (a)--(f), we conclude that the entries in $[\bdv_1,\ldots
,\bdv_{\min(p,q)}]\tran Z$ have the
following properties: all entries below the diagonal are zero; all
entries in the last $q-r$ columns
and on and above the diagonal are independent, where those above the
diagonal are distributed as
standard normal and the square of the $j$th diagonal element
is distributed as $\chi^2_{p-j+1}$.

Notice that if $p\le q$, the $\bdv_j$'s defined above already
constitute a basis for
${\mathbb R}^p$. If $p>q$, we can define $\bdv_j$, $j=q+1,\ldots, p$,
such that they are orthogonal to all columns of $Z$, and to each other.
Define $V_r=[\bdv_{r+1},\ldots, \bdv_p]$. By the nature of eigenvalues,
%
%
\begin{eqnarray}\label{e:natureeig}\qquad
&& \sum_{j=r+1}^{p} \lambda_j(Z Z\tran) \nonumber\\
&&\qquad= \min_{\Phi} \{ \tr(\Phi\tran Z
Z\tran\Phi),\nonumber\\[-8pt]\\[-8pt]
&&\qquad\hspace*{32pt} \mbox{$\Phi$ is a $p\times(p-r)$ matrix with orthonormal columns} \}
\nonumber\\
&&\qquad\le \tr(V_r\tran Z Z\tran V_r)=\sum_{j=r+1}^{p} \bdv_j\tran
Z_2Z_2\tran\bdv_j. \nonumber
\end{eqnarray}
It follows from the summary above that the last expression is a sum of
independent
$\chi^2$ random variable, and a simple calculation shows that the
total degrees of freedom
is $(p-r)(q-r)$.
\end{pf*}
\begin{pf*}{Proof of Theorem \ref{thm:fsir_asymp}}
To study the smallest $m-K_{0}$ eigenvalues of $\widehat V_{(m)}$,
we can equivalently study the smallest squared $m-K_{0}$
singular\vspace*{1pt}
values of $\widehat B_{(m)}$. By the asymptotic theory described in
\citet{r10}
and \citet{r16}, it is
straightforward to show that $\sqrt{n} (\widehat B_{(m)}-B_{(m)})$
converges in distribution. By Theorem 4.1 in \citet{r11}, the
pairwise difference between the smallest $m-K_{(m)}$
singular values\vspace*{2pt} of $\widehat B_{(m)}$ and the singular values
of
$U_n=\CP_2\tran\widehat B_{(m)} \CQ_2$ is $O_p(n^{-3/4})$. So,\vspace*{2pt} for
$K_{(m)}=K_0$, the smallest
$m-K_{0}$ eigenvalues of $\widehat V_{(m)}$ are approximated by the complete
set of
eigenvalues\vspace*{1pt} of $U_nU_n^{\mathrm{T}}$, while, for $K_{(m)}<K_0$, the smallest
$m-K_{0}$ eigenvalues of $\widehat V_{(m)}$ are only approximated by a
subset of eigenvalues of
$U_nU_n^{\mathrm{T}}$. We consider the two cases in more details below.

(i) For $K_{(m)} = K_0$, we will prove (\ref
{eq:eigasymp2})
from which (\ref{eq:eigasymp}) follows easily.
%
%
It follows that $\CQ_2=\CQ_{2\circ}$ and
\[
\CT_{K_0,(m)} = n\tr(U_nU_n\tran)+o_p(1).
\]
%
By (i) of Lemma \ref{lm:approx},
%
%
\begin{equation}\label{eq:asymp_dist}
\CT_{K_0,(m)} \cid\tr(\CZ\Lambda\Xi\Lambda\CZ\tran),
\end{equation}
where $\CZ$ is as given in Lemma \ref{lm:asymp_m} and
$\Xi:=\CJ_\bdg\CQ_2 \CQ_2\tran\CJ_\bdg$.
It is easy to see that $\CJ_\bdg$ and $\CQ_2\CQ_2\tran$ are projection
matrices with rank $H-1$ and $H-K_0$, respectively.
Since $\bdg$ is a column of $\CQ_2$, we have $\CQ_2\CQ_2\tran\bdg
=\bdg$. As a result,
\[
\Xi=\CJ_\bdg\CQ_2 \CQ_2\tran\CJ_\bdg=\CQ_2\CQ_2\tran-\bdg
\bdg\tran,
\]
which is a projection matrix with rank and trace equal to $H-1-K_0$.
Since $\Lambda$ is full rank, $\operatorname{Rank}(\Lambda\Xi\Lambda)=
\operatorname{Rank}(\Xi)=H-K_0-1$. Let $A \Delta A\tran$ be the eigen
decomposition of $\Lambda\Xi\Lambda$ where the column of $A$ are
the orthonormal eigenvectors of $\Lambda\Xi\Lambda$ and
$\Delta=\diag\{\delta_1,\ldots, \delta_{H-K_0-1}\}$ contains the
positive eigenvalues. Write
\[
\tr(Z \Lambda\Xi\Lambda Z\tran)
= \sum_{i=1}^{m-K_0} \bdz_i \Lambda\Xi\Lambda\bdz_i\tran=
\sum_{i=1}^{m-K_0} \bdz_i A\Delta A\tran\bdz_i\tran=
\sum_{i=1}^{m-K_0} \sum_{k=1}^{H-K_0-1} \delta_k\chi^2_{i,k},
\]
where $\bdz_i$ is the $i$th row vector of $Z$, and $\chi_{i,k}^2$
is the $k$th element of $\bdz_iA$. Clearly, the $\chi_{i,k}^2$ are
i.i.d. $\chi^2$ random variables with degree 1.

(ii)
For $K_{(m)} < K_0$, it follows that
\begin{eqnarray*}
\CT_{K_0,(m)}&=& n \times\sum_{j=K_0+1}^m \lambda_j\bigl(\widehat
B_{(m)}\widehat
B_{(m)}\tran\bigr) \\
&=& n \times\sum_{j=K_0-K_{(m)}+1}^{m-K_{(m)}} \lambda_j(U_n
U_n\tran)+o_p(1)\\
&\cid& \sum_{j=K_0-K_{(m)}+1}^{m-K_{(m)}} \lambda_j(ZZ\tran)
\end{eqnarray*}
by Lemma \ref{lm:approx}.
%
Since the last column of $Z$ is identically zero, the last expression
is equal to
\[
\sum_{j=K_0-K_{(m)}+1}^{m-K_{(m)}} \lambda_j\bigl\{ \bigl(Z_*,Z_{\circ
}^{[-1]}\bigr)\bigl(Z_*,Z_{\circ}^{[-1]}\bigr)\tran\bigr),
\]
where $Z_{\circ}^{[-1]}$ denotes the matrix $Z_{\circ}$ minus the
last column. We apply
Proposition~\ref{pp:eign_value_bound}, with $Z_1=Z_*,Z_2=Z_{\circ
}^{[-1]},p=m-K_{(m)},q=H-K_{(m)}-1$
and $r=K_0-K_{(m)}$, to obtain the desired result.
\end{pf*}

\subsection{\texorpdfstring{Proof of Theorem
\protect\ref{thm:fsir_asymp1}}{Proof of Theorem 3.3}}
\begin{lm}\label{lm:columnspace_W}
$W_{(m)}$ has the same column space as $B_{(m)}$.
\end{lm}
\begin{pf}
By definition, $W_{(m)}=B_{(m)} \Lambda(\Lambda\CJ_\bdg
\Lambda)^-$, therefore the column space of $W_{(m)}$ is contained in that
of $B_{(m)}$. Suppose the column rank of $W_{(m)}$ is strictly less
than that of $B_{(m)}$.
Then there exits a nonzero vector $\bdx\in\bbR^m$ such that $\bdx
\tran
W_{(m)}=\mathbf{0}$ but $\bdx\tran B_{(m)}\neq\mathbf{0}$.
Since $B_{(m)}\bdg=\mathbf{0}$, $B_{(m)}\CJ_\bdg=B_{(m)}$ and so
%
%
\begin{equation}\label{e:xTW}
\mathbf{0}=\bdx\tran W_{(m)}= \bdx\tran B_{(m)} \Lambda^{-1}
(\Lambda\CJ_\bdg
\Lambda)(\Lambda\CJ_\bdg\Lambda)^-.
\end{equation}
Observe that $\Lambda^{-1}\bdg$ spans the null space of $\Lambda\CJ
_\bdg\Lambda$.
Since $\bdx\tran B_{(m)}\neq\mathbf{0}$, we conclude that
$\bdx\tran B_{(m)} \Lambda^{-1}= \delta(\Lambda^{-1} \bdg)\tran$
for some constant $\delta$. Thus, $\bdx\tran B_{(m)}=\delta\bdg
\tran$.\vspace*{2pt}
Since $B_{(m)}\bdg=\mathbf{0}$, it follows from (\ref{e:xTW}) that
$\|\bdx\tran B_{(m)}\|^2=\bdx\tran B_{(m)}B_{(m)}\tran\bdx
=\delta\bdg\tran
B_{(m)}\tran\bdx=\mathbf{0}$,
which leads to a contradiction to the assumption that $\bdx\tran
B_{(m)}\neq\mathbf{0}$.
The only possibility left is that the column space of $W_{(m)}$ is the
same as $B_{(m)}$.
\end{pf}
\begin{pf*}{Proof of Theorem \ref{thm:fsir_asymp1}}
As mentioned in the\vspace*{2pt}
proof of Theorem \ref{thm:fsir_asymp}, $ \widehat
B_{(m)}=B_{(m)}+O_p(n^{-1/2})$, which leads to
$\widehat V_{(m)}=V_{(m)}+O_p(n^{-1/2})$ and $\widehat\CP_2 \widehat
\CP{\,}_2\tran=
\CP_2
\CP_2\tran+O_p(n^{-1/2})$. By (\ref{eq:cond_var_est})\vspace*{2pt} and (\ref
{e:vardecomp1}),
we have $\widehat\tau_h=\tau_h +O_p(n^{-1/2})$.
Therefore, $\widehat W_{(m)}$ is a root $n$ consistent estimator of
$W_{(m)}$. The
rest of the proof will follow the same general structure as that of (i) of
Theorem \ref{thm:fsir_asymp}.
Suppose $W_{(m)}$ has the singular-value decomposition
\[
W_{(m)}=\CR
\pmatrix{\widetilde D &0\cr0 &0}
\CS\tran,
\]
where $\CR$ and $\CS$ are, respectively, $m\times m$ and $H\times H$
orthonormal
matrices, and $\widetilde D = \diag( \lambda_1^{1/2}(\Sigma_{(m)}),
\ldots,
\lambda_{K_{(m)}}^{1/2}(\Sigma_{(m)}))$. As before,
consider the partition $\CR=[ \CR_1 | \CR_2 ]$ and $\CS=[ \CS
_1 | \CS_2]$ where
$\CR_{1}$ and $\CS_{1}$ have $K_{(m)}$ columns, and $\CR_2$ and $\CS
_2$ have $m-K_{(m)}$
and $H-K_{(m)}$ columns, respectively. By Lemma
\ref{lm:columnspace_W}, $B_{(m)}$ and $W_{(m)}$ have the same column space,
therefore we can take $\CR_{2}=\CP_{2}$ without loss of generality.
Similar to the definition of $\CQ_2$, we proceed to construct $\CS_2$.
Again, since $\operatorname{span}(\mu_1-\mu,\ldots,\mu_H-\mu)$ has dimension less
than or
equal to $K_0$, there exists a matrix
$\CS_{2\circ}$ with dimension $H\times(H-K_0)$ and orthonormal
columns such that
\[
(\mu_1-\mu,\ldots,\mu_H-\mu)\star(F\Lambda(\Lambda\CJ_\bdg
\Lambda)^-
\CS_{2\circ})=\mathbf{0}.
\]
Observe that
$(\mu_1-\mu,\ldots,\mu_H-\mu)\star F\Lambda(\Lambda\CJ_\bdg
\Lambda)^-\Lambda^{-1}\bdg
= 0$ since $\Lambda^{-1}\bdg$ spans the null space of $\Lambda\CJ
_\bdg\Lambda$.
Without loss of generality, let $\Lambda^{-1}\bdg$ be the last column
of $\CS_{2\circ}$.
Let $T$ be as defined in (\ref{e:T}). As in (\ref{e:BQ2}), we obtain
\begin{eqnarray*}
W_{(m)}\CS_{2\circ} &=& MF\Lambda(\Lambda\CJ_\bdg\Lambda)^-\CS
_{2\circ}\\
&=& T\{(\mu_1-\mu,\ldots,\mu_H-\mu)\} F\Lambda(\Lambda\CJ_\bdg
\Lambda)^-\CS_{2\circ}\\
&=& \mathbf{0}.
\end{eqnarray*}
Since $K_{(m)}=K_0$, we can, and will, take $\CS_2$ to be $\CS
_{2\circ}$. Again, by Theorem 4.1 in \citet{r11}, the
smallest $m-K_{(m)}$ singular values\vspace*{2pt} of $\widehat W_{(m)}$ are
asymptotically equivalent to those of $U_n^\ast:=\CP_2\widehat
W_{(m)}\CS_2$,
so that we have
\[
\CT_{K_0,(m)}^{\ast} = n\tr\{ U_n^\ast(U_n^\ast)\tran\}+o_p(1).
%
\]
Let $\CF=G\CJ_\bdg\Lambda(\Lambda\CJ_\bdg\Lambda)^{-}$,
$\widehat
\CF=\widehat G\CJ_{\widehat\bdg} \widehat\Lambda(\widehat\Lambda
\CJ_{\widehat\bdg}
\widehat\Lambda)^{-}$.
Similar to Lemma \ref{lm:approx},
\[
%
n^{1/2} U_n^\ast=n^{1/2} \CP_2\tran\widehat M \widehat\CF\CS_2= n^{1/2}
\{\CP_2\tran\widetilde M \CF\CS_2+\CP_2\tran(\widehat
M-\widetilde M) \CF\CS
_2\} +o_p(1).
\]
By arguments similar to those in the proof of Lemma \ref{lm:approx},
we have
$n^{1/2} \CP_2\tran(\widehat M-\widetilde M) \CF\CS_2=o_p(1)$. Let
$\Xi^{\ast} = \Lambda\CJ_\bdg\Lambda(\Lambda\CJ_{\bdg}
\Lambda)^{-} \CS_2 \CS_2\tran(\Lambda\CJ_{\bdg} \Lambda
)^{-}\Lambda
\CJ_\bdg\Lambda$.
Thus,
\[
\CT_{K_0,(m)}^{\ast}
=n\tr( \CP_2\tran\widetilde M \CF\CS_2 \CS_2\tran\CF\tran
\widetilde M\tran
\CP_2)+o_p(1)
\cid\tr(Z \Xi^{\ast} Z\tran)
\]
by Lemma \ref{lm:asymp_m}, where $Z$ is $(m-K_0)\times H$ matrix
with independent $\Normal(0,1)$ entries. By the fact that $\Lambda\CJ
_\bdg\Lambda(\Lambda\CJ_{\bdg}\Lambda)^{-}$ is a projection
matrix with only one null vector $\Lambda^{-1}\bdg$ and the
assumption that $\Lambda^{-1}\bdg$ is a column of $S_2$, it easily
follows that $\Xi^{\ast}$
is a projection matrix with trace equal to $H-1-K_0$.
Therefore, $\tr(Z \Xi^{\ast} Z\tran)$ is distributed as $\chi
^2_{(m-K_0)\times(H-K_0-1)}$.
\end{pf*}

\subsection{\texorpdfstring{Proof of Theorem
\protect\ref{thm:increasing_m}}{Proof of Theorem 3.4}}

Some of the variables and matrices introduced in earlier sections
depend on $m$, and we will add the subscript $\mbox{}_{(m)}$ to those
quantities to emphasize this dependence in the proof. Consider the
singular-value decomposition
\[
B_{(m)} = \CP_{(m)}
\pmatrix{D_{(m)} &0\cr0 &0}
\CQ_{(m)}\tran,
\]
where $\CP_{(m)}$ and $\CQ_{(m)}$ have the same partition as before (see
the proof of Theorem~\ref{thm:fsir_asymp}): $\CP_{(m)} = [ \CP
_{1,(m)} |
\CP_{2,(m)} ]$ and $\CQ_{(m)} = [ \CQ_{1,(m)} | \CQ_{2,(m)}
]$. The nonuniqueness of
$\CP_{2,(m)}$ and $\CQ_{2,(m)}$ allows us to construct $\CP_{2,(m)}$
and $\CQ_{2,(m)}$
in a particular way, as follows. It
will be easier to think about the case where $K_{(m)}=K_0$ for $m$
large enough.
We will henceforth make this assumption even though it is not necessary
for the result to hold.
Thus, there exits an ascending sequence $0<m_1<m_2<\cdots
<m_{K_0}<\infty$, such that
\[
m_j=\min\bigl\{m, K_{(m)}\ge j\bigr\},\qquad j=1,\ldots, K_0,
\]
which are the instances where the rank of $B_{(m)}$ changes.

We first construct $\CQ_{2,(m)}$ whose columns span the null row space
of $B_{(m)}$.
Let $\CQ_{2\circ}$ be as in the proof of Theorem \ref{thm:fsir_asymp}.
Define a sequence of orthonormal vectors $\bdq_j, j=1,\ldots,K_0$ by
backward induction:
\begin{eqnarray*}
B_{(m_{K_0}-1)}\bdq_{K_0}&=&\mathbf{0} \quad\mbox{and}\quad
\CQ_{2\circ}^{\mathrm{T}}\bdq_{K_0}=\mathbf{0}; \\
B_{(m_j-1)}\bdq_j&=&\mathbf{0}
\end{eqnarray*}
and
\begin{eqnarray*}
[\bdq
_{j+1},\ldots,
\bdq_{K_0}, \CQ_{2\circ}]^{\mathrm{T}}\bdq_j&=&\mathbf{0},\qquad j=K_0-1,\ldots,2;
\\
{}[\bdq_2,\ldots, \bdq_{K_0}, \CQ_{2\circ}]^{\mathrm{T}} \bdq_1&=&\mathbf{0}.
\end{eqnarray*}
Such a sequence of $\bdq_j$'s clearly exist.
Define
%
%
\begin{equation}\label{eq:Q_2_m_construct}
\CQ_{2,(m)} = \cases{
\bigl[\bdq_{K_{(m)}+1},\ldots, \bdq_{K_0}, \CQ_{2\circ}\bigr], &\quad $m < m_{K_0}$,
\cr
\CQ_{2\circ}, &\quad $m \ge m_{K_0}$.}
\end{equation}
Thus, $\CQ_{2,(m+1)}=\CQ_{2,(m)}$ if $K_{(m+1)}=K_{(m)}$, otherwise
$\CQ_{2,(m+1)}$
equals $\CQ_{2,(m)}$ minus the first column.

We next construct $\CP_{2,(m)}$, a matrix of dimension $m\times(m-K_{(m)})$,
whose columns generate the null column space of $B_{(m)}$.
To do that, we start with $m=K_0+1$ for which we will just make an
arbitrary choice of $\CP_{2,(m)}$ that works.
Suppose we have defined $\CP_{2,(m)}$ for some $m$.
If $K_{(m+1)}=K_{(m)}+1$, let
\[
\CP_{2,(m+1)}= \left[
\matrix{\CP_{2,(m)}\cr{\mathbf0}_{m-K_{(m)}}^{\mathrm{T}}}
\right];
\]
if $K_{(m+1)}=K_{(m)}$, let
%
%
\begin{equation}\label{e:vm+1}\qquad\quad
\CP_{2,(m+1)}=\bigl(\CP_{21,(m+1)}, \bdv_{m+1}\bigr)\qquad \mbox{where }
\CP_{21,(m+1)}= \pmatrix{
\CP_{2,(m)}\cr
\mathbf{0}_{m-K_{(m)}}\tran},
\end{equation}
where $\bdv_{m+1}$ is a new null singular column vector in $\bbR
^{m+1}$. Thus, a whole sequence of $\CP_{2,(m)}$ can
be defined recursively in this manner.

We now briefly summarize some of the key points in the proof of Theorem
\ref{thm:fsir_asymp}.
For each $m\ge K_0+1$, there exists a Gaussian
random matrix $Z_{(m)}$ such that
%
%
\begin{eqnarray}\label{e:joint}
\CP_{2,(m)}\tran\widehat B_{(m)} \CQ_{2,(m)}&\cid& Z_{(m)} \quad\mbox
{and}\nonumber\\[-8pt]\\[-8pt]
\CT_{K_0,(m)}&\cid&\sum_{j=K_0-K_{(m)}+1}^{m-K_{(m)}} \lambda
_j\bigl(Z_{(m)} Z_{(m)}\tran\bigr);\nonumber
\end{eqnarray}
if we write $Z_{(m)}=[ Z_{*,(m)} | Z_{\circ,(m)} ]$ where
$Z_{*,(m)}$ contains the first
$K_0-K_{(m)}$ columns of $Z_{(m)}$, then $Z_{*,(m)}$ is independent of
$Z_{\circ,(m)}$,
and $Z_{\circ,(m)}$ contains independent $\Normal(0,1)$ random
variables except the last
column which contains zeros.
The proof of Proposition \ref{pp:eign_value_bound} shows that there
exist orthonormal vectors
$\bdphi_{1,(m)},\ldots,\bdphi_{m-K_0,(m)}$ in $\bbR^{m-K_{(m)}}$
that are orthogonal to the columns of
$Z_{*,(m)}$ such that
%
%
\begin{equation}\label{e:Xm}
\CX_{(m)}:=\sum_{j=1}^{m-K_{0}} \bdphi_{j,(m)}\tran Z_{(m)}
Z_{(m)}\tran\bdphi_{j,(m)}
\sim\chi^2_{(m-K_0)\times(H-K_0-1)}.
\end{equation}
Note that the $\bdphi_{j,(m)}$'s are obtained by relabeling the $\bdv
_j$'s in that proof. By (\ref{e:joint}) and
(\ref{e:Xm}), using the notion of (\ref{e:natureeig}), we conclude that
$\CT_{K_0,(m)}$ is asymptotically bounded by $\CX_{(m)}$. Thus, we
have the desired stochastic bound
for the first term, $m=K_0+1$, but so far there is nothing new.
To define $\CX_{(m+1)}$, we proceed in a similar
manner by identifying a set of orthonormal vectors $\bdphi_{j,(m+1)},
j=1,\ldots,(m+1)-K_0$.
As we will see, the specific choice of $\CP_{2,(m)}$ and $\CQ
_{2,(m)}$ that was made enables us to directly relate
$Z_{(m)}$ and $Z_{(m+1)}$ in a probability space, and, consequently,
the two bounds as well.

Consider the two situations $K_{(m+1)}=K_{(m)}+1$ and
$K_{(m+1)}=K_{(m)}$ separately.

\textit{Case} 1: $K_{(m+1)}=K_{(m)}+1$.
In view of the relationship between $(\CP_{2,(m)},\break\CQ_{2,(m)})$
and $(\CP_{2,(m+1)}$, $\CQ_{2,(m+1)})$, it is easy to see that
$Z_{(m+1)}$ is equal to $Z_{(m)}$ less the first column.
Below we denote the $k$th column of $Z_{*,(m)}$ by $\bdz_{k,(m)}$.
Define
\[
\bdphi_{j,(m+1)}= \cases{ \bdphi_{j,(m)}, &\quad $j=1,\ldots, m-K_0$, \vspace*{2pt}\cr
\dfrac{(I-\Pi)\bdz_{1,(m)}}{
\|(I-\Pi)\bdz_{1,(m)}\|}, &\quad $j=(m+1)-K_0$,}
\]
where $\Pi$ is the projection matrix onto $\operatorname{span}\{\bdz
_{2,(m)},\ldots, \bdz_{K_0-K_{(m)},(m)}\}$. Observe that the
vectors $\bdphi_{j,(m+1)}, j = 1,\ldots, (m+1)-K_0$, are orthonormal.
Define
%
%
\begin{equation}\label{e:recurbd}\qquad
\CX_{(m+1)}=\sum_{j=1}^{(m+1)-K_{0}} \bdphi_{j,(m+1)}\tran
Z_{(m+1)} Z_{(m+1)}\tran\bdphi_{j,(m+1)}
= \CX_{(m)} + \CX,
\end{equation}
where $\CX=\bdphi_{(m+1)-K_0,(m+1)}\tran Z_{(m+1)} Z_{(m+1)}\tran
\bdphi_{(m+1)-K_0,(m+1)}$.
Note that
\[
\bdphi_{(m+1)-K_0,(m+1)}\in\operatorname{span}^{\perp}\bigl\{\bdz
_{2,(m)},\ldots,\bdz_{K_0-K_{(m)},(m)}\bigr\}
\]
and is independent of $Z_{\circ,(m)}$. The conditioning arguments in
the proof of Proposition \ref{pp:eign_value_bound}
can be applied to conclude that $\CX\sim\chi^2_{H-K_0-1}$ and is
independent of~$\CX_{(m)}$. As a result, $\CT_{K_0,(m)}$ and $\CT
_{K_0,(m+1)}$ are jointly asymptotically bounded by
$\CX_{(m)}$ and $\CX_{(m)} + \CX$.

\textit{Case} 2: $K_{(m+1)}=K_{(m)}$. In this case,
\[
Z_{(m+1)}= \left[
\matrix{Z_{(m)}\vspace*{2pt}\cr\bdw^{\mathrm{T}}}
\right],
\]
where $ \bdw\tran$ is the limit of $\bdv_{m+1}\tran\widehat B_{(m+1)}
\CQ_{2,(m)}$; see (\ref{e:vm+1}).
Define
\[
\bdphi_{j,(m+1)}= \cases{
\left[
\matrix{\bdphi_{j,(m)} \cr0}
\right], &\quad $j=1,\ldots, m-K_0$,
\vspace*{2pt}\cr
\dfrac{(I-\Pi)\bde}{
\|(I-\Pi)\bde\|}, &\quad $j=(m+1)-K_0$,}
\]
where, in this case, $\Pi$ is the projection matrix onto the column
space of
$Z_{*,(m+1)}$ and $\bde=(\mathbf{0}_{m-K_{(m)}}\tran, 1)\tran$.
Define $\CX_{(m+1)}$ as in (\ref{e:recurbd})
and the same jointly asymptotic bound can be concluded for $\CT
_{K_0,(m)}$ and $\CT_{K_0,(m+1)}$, as in the
previous case.

These construction steps can be implemented recursively, thereby
completing the proof of
Theorem \ref{thm:increasing_m}.
\end{appendix}

\section*{Acknowledgments}
We are grateful to the Associate Editor and two referees for their
helpful comments
and suggestions.

\printaddresses

\end{document}